\documentclass[11pt,oneside]{amsart}
\usepackage[a4paper]{geometry}
\usepackage[small,width=\linewidth]{caption}
\usepackage{graphicx}        
\usepackage{overpic}
\usepackage{amssymb}
\usepackage{xypic}

\graphicspath{{./wallner-figures/}}

\newcounter{Satz}

\def\theSatz{\arabic{Satz}}

\newenvironment{Satz}[1]
 {\refstepcounter{Satz}\trivlist\item[\hskip \labelsep{\bf #1\ \theSatz}]\it}
 {\endtrivlist}
\newenvironment{Definition}
 {\refstepcounter{Satz}\trivlist\item[\hskip \labelsep{\bf Definition \theSatz}]\rm}
 {\endtrivlist}
\newenvironment{Example}[1]
 {\refstepcounter{Satz}\trivlist\item[\hskip \labelsep{\bf Example \theSatz.}]
	{\it #1.} \rm}
 {\qed\endtrivlist}

\def\Bullet#1{\Paragraph{$\bullet$~~~#1}}
\def\cput(#1,#2)#3{\put(#1,#2){\hbox to 0pt{\hss{#3}\hss}}}
\def\lput(#1,#2)#3{\put(#1,#2){\hbox to 0pt{\hss{#3}}}}
\def\MM{{\mathord{\mathcal{M}}}}

\def\dMM{\operatorname{\text{\it d}_{\!\MM}}}
\def\TpMM{\mathord{T\hskip-0.4ex_p\hskip-0.1ex\MM}}
\def\wt{\widetilde}
\def\gg{{\mathord{\mathfrak{g}}}}
\def\hh{{\mathord{\mathfrak{h}}}}
\def\ss{{\mathord{\mathfrak{s}}}}
\def\id{{\mathord{\text{\rm id}}}}
\def\ZZ{{\mathbb Z}}
\def\RR{{\mathbb R}}
\def\Lip{{\operatorname{Lip\,}}}
\def\EE{\mathop{{\mathbb E\,}}}
\def\PT{\mathop{{\text{\rm pt}\hskip0.1ex}}\nolimits}
\def\PP{\mathop{{\mathbb P\,}}}
\def\Sp{\mathord{{\text{\it Sp}}}}
\def\Skp{S^{k\!}p}
\def\Spj{{\Sp_{\hskip-0.1ex j}}}
\def\Tp{\mathord{{\text{\it Tp}}}}
\def\GL{\mathord{{\text{\rm GL}}}}
\def\SO{\mathord{{\text{\rm SO}}}}

\def\OO{\mathord{{\text{\rm O}}}}
\def\Pos{\mathord{{\text{\rm Pos}}}}
\def\Sym{\mathord{{\text{\rm Sym}}}}
\def\Paragraph#1{\par\medskip\noindent{\em #1}.}
\def\avg{\operatorname{avg}}

\begin{document}

\title{Geometric subdivision and multiscale transforms}
\author{Johannes Wallner}

\address{Johannes Wallner.  TU Graz, Kopernikusgasse 24, 80180 Graz,
	Austria, email j.wallner@tugraz.at}

\begin{abstract}Any procedure applied to data, and any quantity derived from 
data, is required to respect the nature and symmetries of the data. This
axiom applies to refinement procedures and multiresolution transforms as
well as to more basic operations like averages. This chapter discusses
different kinds of geometric structures like metric spaces, Riemannian
manifolds, and groups, and in what way we can make elementary operations
geometrically meaningful. A nice example of this is the Riemannian metric
naturally associated with the space of positive definite matrices and 
the intrinsic operations on positive definite matrices derived from it. 
We disucss averages first and then proceed to refinement operations
(subdivision) and multiscale transforms. In particular, we report on
the current knowledge as regards convergence and smoothness.
\end{abstract}
\maketitle

\section{Computing averages in nonlinear geometries}

The line of research presented in this chapter was first suggested
by a 2001 presentation by D.\ Donoho on multiscale
representations of discrete data \cite{jw:donoho-lie}.
A subsequent Ph.D.\ thesis and accompanying publication appeared a few
years later \cite{jw:donoho-multi}. Multiscale representations are
intimately connected with refinement procedures (prediction operators).
These are in themselves an interesting topic with applications, e.g.\
in computer graphics. Iterative refinement a.k.a.\ subdivision in turn is
based on the notion of {\em average.}
Consequently this chapter is structured into
the following parts: Firstly a discussion of averages, in particular averages
in metric spaces and in manifolds. Secondly,
subdivision rules and the limits generated by them. Thirdly,
multiresolution representations. 

\index{average|(}
We start with the {\em affine average} w.r.t.\ weights $a_j$ of 
data points $x_j$ contained in a vector space. It is defined by 
	\begin{align}
	x=\avg_{j\in\ZZ} (a_j,x_j) := \sum a_j x_j,\quad \text{where}\quad
	\sum a_j=1.
	\end{align}
In this chapter we stick to finite averages, but we allow negative
coefficients.  For data whose geometry is not that
of a vector space, but that of a surface contained in some Euclidean
space, or that of a group, or that of a Riemannian manifold, this
affine average often does not make sense. In any case it is not natural. 
Examples of such data are, for instance, unit vectors,
positions of a rigid body in space, or the 3 by 3 symmetric positive definite
matrices which occur in diffusion-tensor MRI. In the following paragraphs we show
how to extend the notation of affine average to nonlinear situations
in a systematic way.  We start by pointing out equivalent
characterizations of the affine average:
	\begin{align}
	 x = \avg(a_j,x_j)
	& \iff x \ \text{solves} \  \textstyle \sum  a_j (x_j-x)=0
		\label{jw:eq:implicit}
	\\ &\iff  x = y + \textstyle\sum a_j(x_j-y) \ \text{for any}\ y
		\label{jw:eq:basepoint}
	\\ & \iff x \ \text{minimizes}\ \textstyle \sum a_j \|x-x_j\|^2.
		\label{jw:eq:frechet0}
	\end{align}

\subsubsection*{The Fr\'echet mean}
\index{Frechet mean@Fr\'echet mean}
\index{Riemannian manifold}
\index{Hadamard space}

 Each of \eqref{jw:eq:implicit}--\eqref{jw:eq:frechet0}
has been used to generalize the notion of weighted
average to nonlinear geometries. Some of these generalizations 
are conceptually straightforward. For example,
Equation \eqref{jw:eq:frechet0} has an analogue in any metric
space $(\MM,\dMM)$, namely the weighted {\it Fr\'echet mean} defined
by
	\begin{align}
	\avg_F(a_j,x_j):= \arg\min_{x\in\MM} \sum a_j \dMM(x,x_j)^2.
	\label{jw:eq:frechetmean}
	\end{align}
It is a classical result that in case of nonnegative weights,
the Fr\'echet mean exists and is unique, if
$\MM$ is a Hadamard metric space. This property means $\MM$ is
complete, midpoints exist uniquely, and triangles are slim,
cf.\ \cite{jw:ball-1995}.\footnote{More precisely,
for all $a,b\in\MM $ there is a unique midpoint $x=m(a,b)$ 
defined by $\dMM(x,a)=\dMM(x,b)=\dMM(a,b)/2$,  and 
for any $a,b,c\in\MM$  and points $a',b',c'\in\RR^2$ which have the
same pairwise distances as $a,b,c$, the inequality 
$\dMM(c,m(a,b)) \le d_{\RR^2}(c',m(a',b'))$ holds.}

\Paragraph{The Fr\'echet mean in Riemannian manifolds}
In a surface resp.\ Riemannian manifold $\MM$, the Fr\'echet mean locally
exists uniquely. A main reference here 
is the paper \cite{jw:karcher} by H.~Karcher.  He considered the more
general situation that $\mu$ is a probability measure on $\MM$, where
the mean is defined by
	\begin{align*}
	\avg_F(\mu)= \arg\min_{x\in\MM} \int \dMM(x,\cdot)^2 d\mu.
	\end{align*}
 In this chapter we stick to the elementary case of finite averages
with possibly negative weights. The Fr\'echet mean exists uniquely
if the manifold is Hadamard  -- this property is usually called 
``Cartan-Hadamard'' and is characterized by completeness, simple connectedness,
and nonpositive sectional curvature. For unique existence
of $\avg_F$, we do not even have to require that
weights are nonnegative \cite[Th. 6]{jw:huening-subdiv-2017}.
	\index{Cartan-Hadamard manifold}

\Paragraph{The Fr\'echet mean in the non-unique case}
If the Cartan-Hadamard property is not fulfilled, the Fr\'echet mean
does not have to exist at all, e.g.\ if the manifold is not complete
(cutting a hole in $\MM$ exactly where the mean should
be makes it nonexistent). If the manifold is complete,
the mean exists, but possibly is not unique.

If $\MM$ is complete with nonpositive sectional curvature,
but is not simply connected,
there are situations where a unique Fr\'echet mean of given data points
can still be defined, e.g.\ if the data are connected by a path
$c\colon [a,b]\to\MM$ with $c(t_j)=x_j$. This will be the
case e.g.\ if data represent a time series. Existence or maybe even
canonical existence of such a path depends on the particular application.
We then consider the simply connected
covering $\wt\MM$, find a lifting $\wt c\colon I\to\wt\MM $ of $c$,
compute the Fr\'echet mean $\avg_F(a_j,\wt c(t_j))$, and project it back
to $\MM$. This average does not only depend on the data points and the
weights, but also on the homotopy class of $c$. In fact instead of
a path, any mapping $c\colon I\to\MM$ can be used for such purposes
as long as its domain $I$ is
simply connected \cite{jw:huening-subdiv-2017}.\index{covering, simply connected}

Finally, if $\MM$ is complete but has positive sectional curvatures,
a unique Fr\'echet mean is only defined locally. 
The size of neighbourhoods where uniqueness happens has been 
discussed by 
\cite{jw:dyer2016,jw:dyer2016a,jw:hardering}. This work plays a role in investigating 
convergence of subdivision rules in Riemannian manifolds,
see Section~\ref{jw:ss:subdivmanif}.

\subsubsection*{The exponential mapping}
\index{exponential mapping}

From the different expressions for the affine average, 
\eqref{jw:eq:implicit} and \eqref{jw:eq:basepoint} seem to
be specific to linear spaces, because they involve the $+$ and $-$
operations. However, it turns out that there is a big class of nonlinear
geometries where natural analogues $\oplus$ and $\ominus$ of these
operations exist, namely
the exponential mapping and its inverse. We discuss this construction
in surfaces resp.\ Riemannian manifolds, in groups,
and in symmetric spaces.

\Paragraph{The exponential mapping in Riemannian geometry}
In a Riemannian manifold $\MM$, for any $p\in\MM$ and tangent vector
$v\in \TpMM$, the point $\exp_p(v)$
is the endpoint of the geodesic curve $c(t)$ which starts in $p$,
has initial tangent vector $v$, and whose length equals $\|v\|$.
We let
	\begin{align*}
	p\oplus v &:= \exp_p(v), &  q\ominus p &:=\exp_p^{-1}(q).
	\end{align*}
One property of the exponential mapping is the fact that curves of the form 
	$t\mapsto p \oplus tv$
are shortest paths with initial tangent vector $v$.
The mapping $v\mapsto p\oplus v$ is a diffeomorphism locally around $v=0$.
Its differential equals the identity.

\Paragraph{Properties of the Riemannian exponential mapping}
For complete Riemannian
manifolds, $p\oplus v$ is  always well defined. 
Also $q\ominus p$ exists by the Hopf-Rinow theorem, but
it does not have to be unique. Uniqueness happens if $\dMM(p,q)$ does not
exceed the {\it injectivity radius} $\rho_{\text{inj}}(p)$ of $p$. In 
Cartan-Hadamard manifolds, injectivity radii are infinite and the
exponential mapping does not decrease distances, i.e.,
$\dMM(p\oplus v,p\oplus w)\ge \|v-w\|_{\TpMM}$. The injectivity radius
can be small for topological reasons (e.g.\, a cylinder of small radius
which is intrinsically flat, can have arbitrarily small injectivity
radius), but even in the simply connected
case, one cannot expect $\rho_{\text{inj}}$ to exceed $\pi K^{-1/2}$, 
if $K$ is a positive
upper bound for sectional curvatures.

Further, the $\ominus$ operation and
the Riemannian distance are related by
	\begin{align}
	\label{jw:eq:gradient}
	\nabla \dMM(\cdot,a) (x)  &= -{a\ominus x \over \|a\ominus x\|}, 
	&\nabla \dMM^2(\cdot,a) (x) &= -2(a\ominus x),
	\end{align}
if $v=a\ominus x$ refers to the smallest solution $v$ of $x\oplus a=v$.
For more properties of the exponential
mapping we refer to  \cite{jw:karcher} and to
differential geometry textbooks like \cite{jw:docarmo}.

\Paragraph{The exponential mapping in groups}
In Lie groups, which we describe
only in the case of a matrix group $G$, a canonical
exponential mapping is defined: With the notation
$\gg=T_e G$ for the tangent space in the identity element,
we let 
	\begin{align*}
	v\in \gg\implies
	e\oplus v = \exp(v) = \sum\nolimits_{k\ge 0} {1\over k!} v^k.
	\end{align*}
 The curve $t\mapsto e\oplus tv$ 
is the unique one-parameter subgroup of $G$ whose tangent vector at $t=0$
is the vector $v\in\gg$. Again, $v\mapsto e\oplus v$ is locally a diffeomorphism
whose differential is the identity mapping.

An inverse {\it log} of
{\it exp} is defined locally around $e$.  Transferring the definition
of $\oplus$ to the entire 
group by left translation, the defining relation
$g \oplus gv := g(e\oplus v)$ yields 
	\begin{align*}
	p\oplus v &= p\exp(p^{-1} v), 
	&q\ominus p &= p \log (p^{-1}q ).
	\end{align*}
Addition is always globally well defined, but the difference
$q\ominus p$ might not
exist. For example, in $\GL_n$, the mapping $v\mapsto e\oplus v$ is
not onto. The difference exists always, but not uniquely, in compact groups. 
See e.g.\ \cite{jw:bump-2004-lg}.

\Paragraph{The exponential mapping in symmetric spaces}\index{symmetric space}
Symmetric spaces have the form $G/H$, where 
$H$ is a Lie subgroup of $G$. There are several definitions which are
not entirely equivalent. We use the one that 
the tangent spaces $\gg=T_e G$, $\hh=T_e H$
obey the condition that $\hh$ is the $+1$
eigenspace of an involutive Lie algebra automorphism $\sigma$ of 
$\gg$.\footnote{i.e., $\sigma$ obeys the law
$\sigma([v,w])=[\sigma(v),\sigma(w)]$, where in the matrix group case,
the Lie bracket operation is given by $[v,w]=vw-wv$.}
The tangent space $\gg/\hh$ of $G/H$ in the point $eH\in G/H$ is naturally
identified with the $-1$ eigenspace $\ss$ of the involution, and is
transported to all points of $G/H$ by left translation. The exponential
mapping in $G$ is projected onto $G/H$
in the canonical way and yields the exponential mapping in the symmetric
space.

We do not go into more
details but refer to the comprehensive classic
\cite{jw:helgason} instead.
Many examples of well-known manifolds fall into this category, e.g.\
the sphere $S^n$, hyperbolic space $H^n$, and the Grassmannians. We give
an important example:

\begin{Example}{The Riemannian symmetric space of positive-definite matrices}%
\label{jw:ex:pos}
The space $\Pos_n$ of positive definite 
$n\times n$ matrices is made a metric space by letting 
	\begin{align}
	d(a,b) = 
		\|
			\log(a^{-1/2} b a^{-1/2})
		\|_2
	= 
	\Big(\sum\nolimits_{\lambda_1,\ldots,
		\lambda_n\in\sigma(a^{-1}b)} \log^2\lambda_j\Big)^{1/2}.
	\label{jw:metric:Pos_n} 
	\end{align}
 Here $\|\cdot\|_2$ means the
Frobenius norm, and $\sigma(m)$ means the eigenvalues of a matrix.

The metric \eqref{jw:metric:Pos_n} 
is actually that of a Riemannian manifold.  $\Pos_n$, as an open
subset of the set $\Sym_n$ of symmetric matrices, in each point has
a tangent space $T_a\Pos_n$ canonically isomorphic to $\Sym_n$ as a
linear space.  The Riemannian metric in this space is defined by
$\|v\| = \|a^{-1/2} va^{-1/2}\|_2$. 

$\Pos_n$ is also a symmetric space:
We know that any $g\in\GL_n$ can be
uniquely written as a product $g=au$, with $a=\sqrt{gg^T}\in\Pos_n$
and $u\in\OO_n$.  Thus $\Pos_n=G/H$, with 
$G=\GL_n$, $H=\OO_n$, and the canonical projection $\pi(x)=\sqrt{xx^T}$.

The respective tangent spaces $\gg,\hh$ of $G,H$ are 
given by $\gg=\RR^{n\times n}$ and $\hh={\mathfrak{so}}_n$, which
is the set of skew-symmetric $n\times n$ matrices.
The involution $\sigma(x)=-x^T$ in $\gg$ obeys
$[\sigma(v),\sigma(w)]=\sigma([v,w])$, and $\hh$ is its $+1$
eigenspace. We have thus recognized $\Pos_n$ as a symmetric space.
It turns out that 
$a\oplus v =a\exp(a^{-1}v)$, where {\it exp} is the
matrix exponential function.

The previous paragraphs define two different structures on
$\Pos_n$, namely that of a Riemannian manifold, and that of a symmetric 
space. They are 
compatible in the sense that the $\oplus$, $\ominus$ operations derived
from either structure coincide. For more information we refer to
\cite{jw:moakher,jw:fletcher-2004,jw:welk-2006}.
Subdivision in particular is treated
by \cite{jw:itai-sharon}.
\end{Example}

\subsubsection*{Averages defined in terms of the exponential mapping}

If $\oplus$ and $\ominus$ are defined as discussed in the previous
paragraphs, it is possible to define a weighted affine average implicitly
by requiring that 
	\begin{align}
	x=\avg_E(a_j,x_j) :\iff \sum a_j (x_j\ominus x)=0.
	\label{jw:eq:implicitaverage}
	\end{align}
Any Fr\'echet mean in a Riemannian manifold is also an average in this sense,
which follows directly from \eqref{jw:eq:frechetmean} together with
\eqref{jw:eq:gradient}.  Locally, $\avg_E$ is well defined and unique.
As to the size of neighbourhoods where this happens, in the Riemannian
case the proof given by \cite{jw:dyer2016,jw:dyer2016a} for
certain neighbourhoods
enjoying unique existence of $\avg_F$ shows that the very same
neighbourhoods also enjoy unique existence of $\avg_E$.

\Paragraph{Affine averages with respect to a base point}\index{base point}
From the different expressions originally given for the affine
average, $x=y+\sum a_j(x_j-y)$ is one we have not yet defined
a manifold analogue for.  With $\ominus$ and $\oplus$ at our disposal,
this can be done by
	\begin{align}
	\avg_y(a_j;x_j) :=
	y\oplus\sum a_j(x_j\ominus y).
	\label{jw:eq:basepoint:manif}
	\index{log/exp subdivision}
	\end{align}
We call this the log/exp average with respect to the base point $y$.
It has the disadvantage of a dependence on the base point, but for
the applications we have in mind, 
there frequently is a natural choice of base point. Its advantages lie in
the easier analysis compared to the Fr\'echet mean. One should also
appreciate that the Fr\'echet mean is a log/exp mean w.r.t.\ to a basepoint,
if that basepoint is the Fr\'echet mean itself:
	\begin{align}
	y=\avg_F(a_j;x_j) \implies \avg_y(a_j;x_j)=
		y \oplus \sum a_j(x_j\ominus y) = y\oplus 0 = y,
	\label{jw:eq:self}
	\end{align}
because of \eqref{jw:eq:implicitaverage}.
This may be a trivial point, but it has been essential in proving smoothness
of limit curves for manifold-based subdivision processes
(see Th.\ \ref{jw:th:philipp} and \cite{jw:grohs-2010-gpa}).

The possibility to define averages w.r.t.\ basepoints rests on the 
possibility of defining $\ominus$, which has been discussed above.

\index{average|)}

\section{Subdivision}
\label{jw:sec:subdivision}
\index{subdivision|(}

\subsection{Defining stationary subdivision}

Subdivision is a refinement process acting on input data lying in some
set $\MM$, which in the simplest case are indexed over
the integers and are interpreted as samples of a function $f\colon
\RR\to\MM$.  A subdivision
rule {\em refines} the input data, producing 
a sequence $\Sp$ which is thought of denser samples of either $f$ itself, or
of a function approximating $f$.

	\index{dilation factor}
One mostly considers binary rules, whose application 
``doubles'' the number of data points. The {\em dilation factor} of the rule,
generally denoted by the letter $N$, then equals $2$.  
We require that the subdivision rule is invariant w.r.t.\ index shift,
which by means of the left shift operator $L$ can be formalized as
	\begin{align*}
	L^N S = S L.
	\end{align*}
We require that each point $\Sp_i$ depends only on finitely many
data points $p_j$. Together with shift invariance this means that
there is $s>0$ such that $p_i$ influences only
$\Sp_{Ni-s},\ldots,\Sp_{Ni+s}$.

Subdivision rules are to be iterated: We create finer and finer data
	\begin{align*}
	p, \ \Sp, \ S^2 p,\ S^3 p,\ \ldots,
	\end{align*}
 which we hope approach a continuous limit (the proper definition of
which is given below).

	\index{corner cutting}
 Subdivision was invented by G.\ de Rham \cite{jw:derham-1962}, who
considered the process of iteratively cutting corners from a convex
polygon contained in $\MM=\RR^2$,
and asked for the limit shape. If cutting corners is done
by replacing each edge $p_ip_{i+1}$ by the shorter edge with
vertices $\Sp_{2i}=(1-t)p_i + tp_{i+1}$, $\Sp_{2i+1}=tp_i+(1-t)p_{i+1}$, 
this amounts to a subdivision rule.  
In de Rham's example, only two data points $p_i$ contribute to any 
individual $\Spj$. 

	\index{primal subdivision}
	\index{dual subdivision}
\Paragraph{Primal and dual subdivision rules}
The corner-cutting rules mentioned above
are invariant w.r.t.\ reordering indices
according to $\ldots, 0\mapsto 1, 1\mapsto 0, 2\mapsto -1,\ldots$. With
inversion $U$ defined by $(Up)_i=p_{-i}$  we can write this invariance
as $(LU)S=S(LU)$. An even simpler kind of symmetry is enjoyed  by
subdivision rules with obey $US=SU$. The latter are called primal rules,
the former dual ones.  The reason why
we emphasize these properties is that they give guidance
for finding manifold analogues of linear subdivision rules.

\Paragraph{Subdivision of multivariate data}
It is not difficult to generalize the concept of subdivision to 
multivariate data $p\colon \ZZ^s\to\MM$ indexed over the 
grid $\ZZ^s$. A subdivision rule $S$ must fulfill
$L_v^N S = SL_v$, for all shifts $L_v$ w.r.t.\ a vector $v\in\ZZ^s$.

Data with combinatorial singularities have to be treated separately,
cf.\ Sec.~\ref{jw:sec:irregular}.
Here basically only the bivariate case is studied, but this has been
done extensively, mostly because of applications in Computer Graphics
\cite{jw:Peters2008}.

\subsubsection*{Linear subdivision rules and their nonlinear analogues}

A linear subdivision rule acting on data $p\colon \ZZ^2\to\RR^d$
has the form
	\begin{align*}
	\Sp_i = \sum\nolimits_{j} a_{i-Nj} p_j.
	\end{align*}
If the sum $\sum_j a_{i-Nj}$ of coefficients
contributing to $\Sp_i$ equals $1$, the application of
the rule amounts to computing a weighted average:
	\begin{align}
	\Sp_i = \avg(a_{i-Nj};p_j).
	\label{jw:eq:linearrule}
	\end{align}
 Subdivision rules not expressible in this way might occur
as auxiliary tools in proofs, but are not 
meant to be applied to data which are {\it points} of an affine space.
This is because if
$\sum a_{i-Nj}\ne 1$, then the linear combination $\sum a_{i-Nj}p_j$
is not translation-invariant, and the rule depends on the choice
of origin of the coordinate system.

Besides, the iterated application of rules not expressible as weighted
averages either leads to divergent data $\Skp$, or alternatively, to 
data approaching zero.
For this reason, one exclusively considers linear rules of the form
\eqref{jw:eq:linearrule}. A common definition of {\em convergent}
subdivision rule discounts the case of zero limits and recognizes
translation invariance as a necessary condition for convergence, cf.\
\cite{jw:dyn:1992:SSIC}.

For a thorough treatment of linear subdivision rules,  conveniently
done via $S$ acting as a linear operator in $\ell^\infty(\ZZ^s,\RR)$ and using
the appropriate tools of approximation theory, see, e.g.\
\cite{jw:cavaretta-1991-ss}. 

In the following we discuss some nonlinear, geometric, versions of
subdivision rules. We use the various nonlinear versions of averages
introduced above, starting with the Fr\'echet mean in metric spaces.

\Bullet{Subdivision using the F\'echet mean}
A natural analogue of \eqref{jw:eq:linearrule} is found by
replacing the affine average by the Fr\'echet mean. This 
procedure is particularly suited for Hadamard metric spaces and also
in complete Riemannian manifolds. 

\Bullet{Log/exp subdivision}
In a manifold equipped with an exponential mapping, an
	analogue of \eqref{jw:eq:linearrule} is defined by
	\begin{align*}
	\Tp_i = \avg_{m_i}(a_{i-Nj};p_j),
	\end{align*}
	where $m_i$ is a base point computed in a meaningful manner from
	the input data, e.g.\ $m_i=p_{\lfloor i/N\rfloor}$.
	In case of combinatorial 
	symmetries of the subdivision rule, it makes sense to make the
	choice of $m_i$ conform to these symmetries.

\Bullet{Subdivision using projections}
If $\MM$ is a surface embedded in a vector space and 
$\pi$ is a projection onto $\MM$, we might use the subdivision rule
	\begin{align*}
	\Tp_i = \pi(\Sp_i).
	\end{align*}
	If the intrinsic symmetries of $\MM$ extend to symmetries of
	ambient space, then this {\em projection analogue}
	of a linear subdivision rule is  even intrinsic -- see
	Example~\ref{jw:ex:motion}.

\begin{figure}[t]
\begin{overpic}[width=0.30\textwidth]{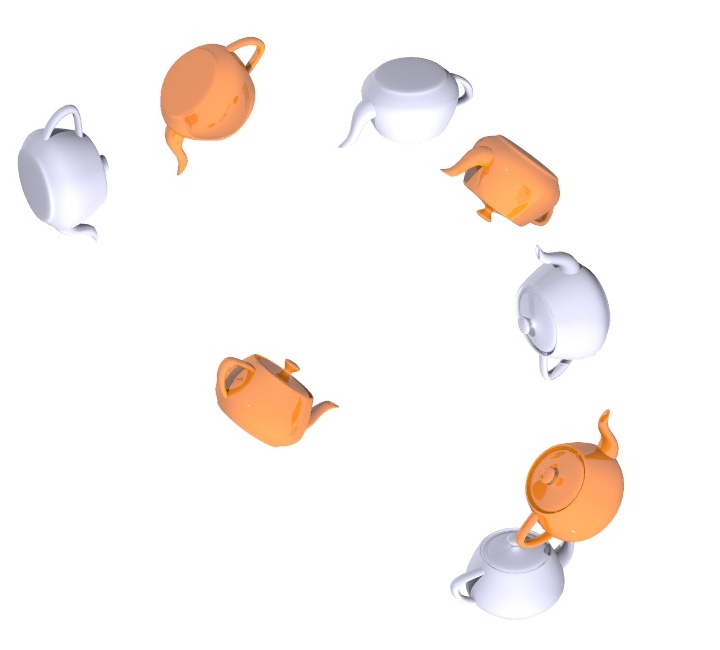}
	\cput(9,63){$p_0=p_4$}
	\cput(56,73){$p_1$}
	\cput(79,44){$p_2$}
	\cput(72,9){$p_3$}
	\cput(27,88){$\Tp_1$}
	\put(77,70){$\Tp_3$}
	\put(90,22){$\Tp_5$}
	\lput(30,30){$\Tp_7$}
	\cput(50,50){\fbox{$p,\Tp$}}
\end{overpic}\relax
\begin{overpic}[width=0.30\textwidth]{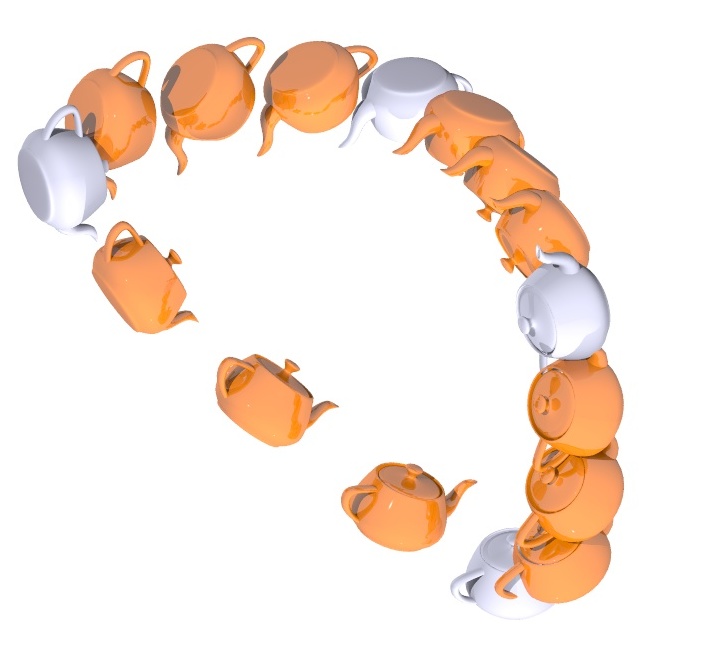}
	\cput(50,50){\fbox{$p,T^2p$}}
\end{overpic}\hfill
\begin{minipage}[b]{.37\textwidth}
	\scriptsize\begin{align*}
 	\Sp_{2i}&=p_i
	\\
	\Sp_{2i+1}&={9\over 16}(p_i+p_{i+1})
		-{1\over 16}(p_{i-1}+p_{i+2})
	\\
	T &= \pi \circ S
	\\
	\end{align*}
\end{minipage}
\caption{Subdivision by projection in the motion group
$\RR^3\rtimes\OO_3$. A 4-periodic sequence $p_i=(c_i,u_i)$ of positions of a rigid
body is defined by the center of mass $c_i$, and an
orientation $u_i\in\OO_3$. Both components undergo subdivision w.r.t.\ the
interpolatory four-point rule $S$,
where the matrix part is subsequently projected
back onto $\OO_3$ in an invariant manner.}
\label{jw:fig:teapot}
\end{figure}

\begin{Example}{Subdivision in the motion group}\label{jw:ex:motion}\relax
The groups $\OO_n$  and $\SO_n$ are ${1\over 2}n(n-1)$-dimensional
surfaces in the linear space
$\RR^{n\times n}$.  A projection onto
$\OO_n$ is furnished by singular value decomposition, or in an alternate
way of expressing it, by the 
polar decomposition of Example~\ref{jw:ex:pos}:
	\begin{align*}	
	\pi\colon \GL_n\to \OO_n, \
	\pi(g)= (gg^T)^{-1/2} g.
	\end{align*}
This projection is $\OO_n$-equivariant in the sense that for $u\in\OO_n$,
we have both
$\pi(ug)=u\pi(g)$ and $\pi(gu)=\pi(g)u$. The same invariance applies
to application of a linear
subdivision rule acting in $\RR^{n\times n}$.
So for any given data in $\OO_n$, and a linear subdivision rule $S$,
the subdivision rule $\pi\circ S$ produces data in $\OO_n$ in a 
geometrically meaningful way,  as long as we do not exceed the bounds
of $\GL_n$. Since $\GL_n$ is a rather big neighbourhood of $\OO_n$, this
is in practice no restriction. Figure \ref{jw:fig:teapot}
shows an example.
\end{Example}

\subsection{Convergence of subdivision processes}

\Paragraph{Definition of convergence}
When discrete data $p$ are interpreted as samples of 
a function, then refined data $\Sp$, $S^2 p$ etc.\ are interpreted
as the result of sampling which is $N$ times, $N^2$ times etc.\ as dense
as the original.  We therefore define a convergent refinement
rule as follows.

\begin{Definition}\label{jw:defn:convergence}\relax
Discrete data $\Skp\colon \ZZ^s\to\MM$ at
the $k$-th iteration of refinement
determine  a function $f_k\colon N^{-k}\ZZ^s\to\MM$,
whose values are the given
data points: For any $N$-adic point $\xi$, we have
$(\Skp)_{N^k \xi} = f_k(\xi)$, provided $N^k \xi$ is an integer. For all
such $\xi$, the sequence $(f_k(\xi))_{k\ge 0}$ is eventually
defined and we let $f(\xi)=\lim_{k\to\infty} f_k(\xi)$. 
We say $S$ is convergent for input data $p$, if 
the limit function $f$ exists for all $\xi$ and is 
continuous. It can be uniquely extended to a continuous
function $S^\infty p\colon \RR^s\to\MM$.
\end{Definition}

Another way of defining the limit is possible if data
$p_i, \Sp_i,\ldots$ lie in a vector space. We linearly
interpolate them by functions $g_0,g_1,\ldots$ with
$g_k(N^{-k}i)=S^k p_i$. Then the limit 
of functions $g_k$ agrees with the limit of
Def.~\ref{jw:defn:convergence} (which is pointwise, 
but in fact convergence is usually uniform on compact sets.)

The following lemma is the basis for investigating convergence of
subdivision rules in metric spaces. The terminology is that of 
\cite{jw:dyn-sharon17a,jw:dyn-sharon17b}.

\begin{Satz}{Lemma} \label{jw:lemma1}
Let $\MM$ be a complete metric space,  and let the
subdivision rule $S$ operate with dilation $N$
on data $p\colon \ZZ^s\to \MM$.
We measure the density of the data by
	\begin{align*}
	\delta(p) = \sup_{|i-j|\le 1} \dMM(p_i,p_j),
	\end{align*}
where we use the 1-norm on the indices.
$S$ is contractive, resp.\ displacement-safe, if 
	\begin{align*}
	 \delta(\Sp)\le \gamma \delta(p),  
			\ \ \text{for some}\ \gamma<1,
	\ \ \text{resp.} \ \
	 \sup\nolimits_{i\in\ZZ^s}\dMM(\Sp_{Ni},p_i)\le \lambda \delta(p). 
	\end{align*}
If these two conditions are met,
any input data with bounded density
have a limit $S^\infty p$,
which is H\"older continuous with exponent $-{\log\gamma\over \log N}$.
\end{Satz}

\begin{proof} Contractivity implies $\delta(\Skp)\le\gamma^k\delta(p)$.
For any $N$-adic rational point $\xi\in N^{-r}\ZZ^s$, the sequence
$f_k(\xi) = (\Skp)_{N^k\xi}$ is defined for all $k\ge r$. It is 
Cauchy, since 
	\begin{align*}
	\dMM(f_k(\xi), f_{k+1}(\xi))
	\le
	\lambda \delta(\Skp) \le \lambda\gamma^k\delta(p).
	\end{align*}
Thus the limit function $S^\infty p \equiv f$
is defined for all $N$-adic points.

Consider now two $N$-adic points $\xi,\eta$. Choose
$k$ such that $N^{-(k+1)}\le|\xi-\eta|\le N^{-k}$. For all $j\ge k$, 
approximate $\xi$ resp.\ $\eta$ by $N$-adic points $a_j,b_j\in N^{-j}\ZZ^s$,
such that none of $|a_j-a|$, $|b_j-b|$, $|a_j-a_{j+1}|$ $|b_j-b_{j+1}|$
exceeds $sN^{-j}$. One can choose $a_k=b_k$. The sequence
$a_j$ is eventually constant with limit $\xi$, and similarly the
sequence $b_j$ is eventually constant with limit $\eta$. 
Using the symbol $(*)$ for ``similar terms involving $b_j$ instead
of $a_j$'', we estimate 
\def\LHS{
	\dMM(f(\xi),f(\eta))
}
	 \begin{align*}
	\LHS
	&\le
	\sum\nolimits_{j\ge k}\dMM(f_j(a_j),f_{j+1}(a_{j+1}))
	+(*)
	\\&\le
		\sum
		\dMM(f_j(a_j),f_{j+1}(a_{j}))+
		\dMM(f_{j+1})(a_j),f_{j+1}(a_{j+1})) + (*) 
	.
	\end{align*}
 Using the contractivity and displacement-safe condition, we further
get
	\begin{align*}
	\LHS &\le
	 2 \sum\nolimits_{j\ge k} \lambda\delta(S^jp
		)+ s\delta(S^{j+1}p) 
	\\[-1ex]&\le
	2(\lambda+s\gamma)\delta(p)\sum\nolimits_{j\ge k} \gamma^j 
	\le C \delta(p) {\gamma^k\over 1-\gamma}.
	\end{align*}
The index $k$ was chosen such that $k \le  -\log|\xi-\eta| /\log N$, so
in particular $\gamma^k \le \gamma^{-\log|\xi-\eta|/\log N}$. We 
conclude that
	\begin{align*}
	\dMM(f(\xi),f(\eta))
	\le
	C'  \gamma^{-\log|\xi-\eta|/\log N}
	=
	C' |\xi-\eta|^{-\log\gamma/\log N}.
	\end{align*}
Thus $f$ is continuous with H\"older exponent
$-{\log\gamma\over \log N}$ on the $N$-adic rationals, and so is
the extension of $f$ to all of $\RR^s$.
\end{proof}

\noindent
The scope of this lemma can be much expanded by some obvious modifications.

\Bullet{Input data with unbounded density $d(p)$} Since
points $\Spj $ only depend on finitely many $p_i$'s, 
there is $m>0$ such that $p_i$ only influences $\Sp_{Ni+j}$ with
$|j|<m$. By iteration, $p_i$ influences $S^2p_{N^2i+j}$ with
$|j|<Nm+m$, and so on. It follows that $p_i$ influences the value 
$S^\infty p(i+\xi)$ of the limit function only for $|\xi|<{m\over N}+
{m\over N^2}+\cdots = {m\over N-1}$. We can therefore easily analyze the
restriction of the limit function to some box by re-defining all
data points away from that box in a manner which makes $d(p)$ finite.

\Bullet{Partially defined input data}
If data are defined not in all of 
$\ZZ^s$ but only in a subset, the limit function is defined for a certain
subset of $\RR^s$. Finding this subset goes along the same lines as the
previous paragraph -- we omit the details.

\Bullet{Convergence for special input data}
In order to check convergence for particular input data $p$, it is
sufficient that the contractivity and displacement-safe conditions
of Lemma~\ref{jw:lemma1}
hold for all data $\Skp$ constructed by iterative
refinement from $p$. A typical instance of this
case is that contractivity can be shown only if $\delta(p)$ does not
exceed a certain threshold $\delta_0$.
It follows that neither does $\delta(\Skp)$,
	and Lemma~\ref{jw:lemma1} applies to all $p$ with $\delta(p)\le\delta_0$.

\Bullet{Powers of subdivision rules}
A subdivision rule $S$ might enjoy convergence like a 
contractive rule without being contractive itself. This phenomenon
is analogous to a linear operator $A$ having norm $\|A\|\ge 1$ but spectral
radius $\rho(A)<1$, in which case some $\|A^m\|<1$.
In that case we consider some power
$S^m$ as a new subdivision rule with dilation factor $N^m$.
If $S^m$ is contractive with factor $\gamma^m<1$,
Lemma~\ref{jw:lemma1} still applies,  and limits enjoy H\"older smoothness with
exponent $-{\log\gamma^m\over \log N^m} = -{\log\gamma\over \log N}$.

\begin{Example}{Convergence of linear subdivision rules}\relax
\label{jw:ex:linear}\relax
Consider a  univariate subdivision rule $S$ defined by
finitely many nonzero coefficients $a_j$ via \eqref{jw:eq:linearrule}.
$S$ acts as a linear operator on sequences
$p\colon\ZZ\to \RR^d$. The norm
$\|p\|=\sup_i\|p_i\|_{\RR^d}$ induces an operator norm $\|S\|$ which obeys 
$\|\Sp\| \le \|S\| \|p\|$.  It is an exercise to check
$\|S\|=\max_i\sum_j|a_{i-Nj}|$.  Equality is attained
for suitable input data with values in $\{-1,0,1\}$.

With $(\Delta p)_i=p_{i+1}-p_i$ we express the density of the data as
$\delta(p)=\sup\|\Delta p_i\|$. 
Contractivity means  that
$\sup \|\Delta \Sp_i\| \le \gamma \sup \|\Delta p_i\|$ for some $\gamma<1$.

	\index{derived subdivision rule}
Analysis of this contractivity condition uses a trick based on the generating functions
$p(z)=\sum p_j z^j$  and
$a(z)=\sum a_j z^j$. Equation
\eqref{jw:eq:linearrule} translates to the relation
 $(\Sp)(z)=a(z)p(z^N)$ between generating functions, and we also have 
$\Delta p(z)=(z^{-1}-1)p(z)$.  The trick consists in introducing 
the {\em derived} subdivision rule $S^*$ with coefficients $a_j^*$ which obeys
$S^*\Delta=N\Delta S$. The corresponding relation between generating functions
reads
	\begin{align*}a^*(z)\Delta p(z^N)
		& =N(z^{-1}-1)a(z)p(z^N)
		\iff a^*(z)(z^{-N}-1) = N(z^{-1}-1)a(z)
	\\
	\iff a^*(z) &= Na(z)z^{N-1}{z-1\over z^N-1} 
		= Nz^{N-1}{a(z)\over 1+z+\cdots+z^{N-1}}.
	\end{align*}
This division is possible in the ring of Laurent polynomials, because  for
all $i$,
$\sum_j a_{i-Nj}=1$. The contractivity condition now reads $\sup\|\Delta \Sp_i\|
={1\over N} \sup\|S^*\Delta p_i\|\le {1\over N}\|S^*\| \sup\|\Delta p_i\|$,
i.e., the contractivity factor of the subdivision rule $S$ is bounded from
above by ${1\over N}\|S^*\|$.  
The ``displacement-safe'' condition of 
Lemma \ref{jw:lemma1} is fulfilled also, which we leave as an exercise
(averages of points $p_i$ are not far from the $p_i$'s).

The above computation leads to a systematic procedure for checking convergence: we 
compute potential contractivity
factors ${1\over N}\|S^*\|$, ${1\over N^2}\|S^{2*}\|$, and so on, until
one of them is $<1$.  
The multivariate case is analogous but more complicated 
\cite{jw:dyn-2002-ss,jw:dyn:1992:SSIC,jw:cavaretta-1991-ss}.
\end{Example}

\begin{Example}{Convergence of geodesic corner-cutting rules}\relax
\label{jw:ex:cornercutting}\relax
Two points $a,b$ of a complete Riemannian manifold $\MM$ are joined
by a shortest geodesic path $t\mapsto a\oplus tv$, $v=b\ominus a$,
$t\in[0,1]$. The difference vector $v$ and thus the path
are generically unique, but do not have to be, if the distance between
$a$ and $b$ exceeds both injectivity radii $\rho_{\text{inj}}(a),
\rho_{\text{inj}}(b)$. 
The point $x=a\oplus tv$ has 
$\dMM(a,x)=t\dMM(a,b)$,
$\dMM(b,x)=(1-t)\dMM(a,b)$. It 
is a Fr\'echet mean of points $a,b$ w.r.t.\ weights $(1-t),t$. 

	\index{corner cutting}
With these preparations, we consider two elementary operations
on sequences, namely {\em averaging} $A_t$ and {\em corner cutting} $S_{t,s}$:
	\begin{align*}
	(A_t p)_i &= 
			p_i \oplus t (p_{i+1}\ominus p_i),
	&(S_{ts}p)_j &= \begin{cases}
		p_i \oplus t (p_{i+1}\ominus p_i)
			& \text{if}\ j=2i, \\
		p_i \oplus s (p_{i+1}\ominus p_i)
			& \text{if}\ j=2i+1.
		\end{cases}
	\end{align*}
The distance of $A_tp_i$ from $A_t p_{i+1}$ is bounded by the 
length of the broken geodesic path which connects the first point with
$p_{i+1}$ and continues on to the second; its length is bounded by 
$\delta(p)$.
Similarly, the distance of successive points of the sequence 
$S_{ts}p$, for $0\le t<s\le 1$ is estimated by $\max(1-(s-t),s-t)\delta(p)$.
It follows immediately that a concatenation of operations of this kind
is a subdivision rule where Lemma \ref{jw:lemma1} applies,
if at least one $S_{t,s}$ with $0<s-t<1$ is involved. Any such
concatenation therefore is a convergent subdivision rule  in any
complete Riemannian manifold.
A classical example are the rules $S^{(k)}=(A_{1/2})^k\circ S_{0,1/2}$,
which insert midpoints,
	\begin{align*} \textstyle
	S^{(1)}p_{2i} &= p_i, 
	&S^{(1)}p_{2i+1} &= p_i\oplus{1\over 2}(p_{i+1}\ominus p_i),
	\end{align*}
and then compute $k$ rounds of averages. E.g.,
	\begin{align*}
	\begin{array}{*9l}
	S^{(2)}p_{2i}& =S^{(1)}p_{2i}
		& \oplus & {1\over 2}(S^{(1)}p_{2i+1}
		&\ominus & S^{(1)}p_{2i})
	&=p_i\oplus {1\over 4}(p_{i+1}\ominus p_i),
	\\[1.5ex]
	S^{(2)}p_{2i+1}& =S^{(1)}p_{2i+2}
		& \oplus & {1\over 2}(S^{(1)}p_{2i+2}
		&\ominus & S^{(1)}p_{2i+1})
	&=p_i\oplus {3\over 4}(p_{i+1}\ominus p_i).
	\end{array}
	\end{align*}
The rule $S^{(2)}$ (Chaikin's rule, see \cite{jw:chaikin}) is one of de Rham's
corner cutting rules.
In the linear case, $S^{(k)}$ has
coefficients $a_j={1\over 2^{k}} {k\choose j}$, apart from an index shift.
Its limit curves are the B-spline curves whose control points are the
initial data $p_j$ \cite{jw:riesenfeld1975}.

The corner-cutting rules discussed above are well defined and convergent in any
Hadamard metric space -- those spaces have geodesics in much the
same way as Riemannian manifolds. Subdivision rules based on
geodesic averaging (not necessarily restricted to values
$t,s\in[0,1]$) have been treated by
\cite{jw:wallner-2004-cca,jw:wallner-2006-is,jw:dyn-sharon17a,jw:dyn-sharon17b}.
We should also mention that adding a round $A_{1/2}$ to a subdivision
increases smoothness of limit curves, which was recently confirmed
in the manifold case \cite{jw:yu-smoothing-2018}.
\end{Example}

\begin{Example}{Convergence of interpolatory rules}\index{interpolatory subdivision}\label{jw:ex:conv:interpolatory}\relax
\label{jw:ex:interpolatory}\relax
A subdivision rule $S$ with dilation factor $N$
is called {\it interpolatory} if $\Sp_{Ni}=p_i$, i.e., the old
data points are kept and new data points are inserted in between. In the
linear case, a very well studied subdivision rule of this kind is
the four-point rule proposed by Dyn, Gregory and Levin \cite{jw:dyn:1987:4p}.
We let $\Sp_{2i}=p_i$ and 
	\begin{align*}
	\textstyle
	\Sp_{2i+1} 
	&= -\omega p_{i-1} \textstyle
		 + ({1\over 2}+\omega)p_i
		 + ({1\over 2}+\omega)p_{i+1}
		-\omega p_{i+2}
	\\
	&= {p_i+p_{i+1}\over 2}
		-\omega \Big(p_{i-1} - {p_i+p_{i+1}\over 2}\Big)
		-\omega \Big(p_{i+2} - {p_i+p_{i+1}\over 2}\Big)
	.
	\end{align*}
In the special
case $\omega={1\over 16}$, the point $\Sp_{2i+1}$ is found by evaluating
the cubic Lagrange polynomial interpolating $p_{i-1},\ldots, p_{i+2}$, which accounts 
for the high approximation order of $S$. There is in fact a whole
series of interpolatory rules based on the idea of evaluating 
Lagrange interpolation polynomials (the Dubuc-Deslauriers subdivision
schemes, see \cite{jw:deslauriers:1989:siip}). 

\begin{figure}[t]
\begin{overpic}[width=0.48\textwidth]{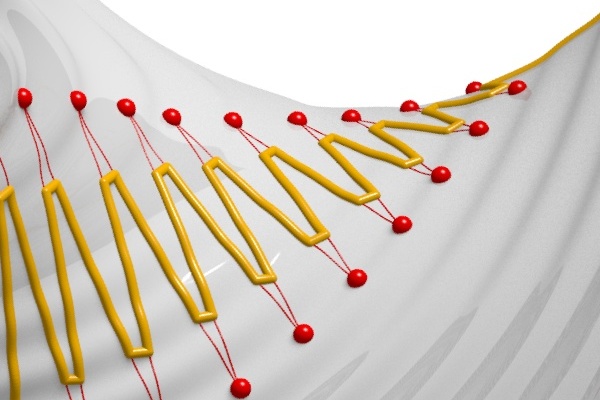}
	\lput(98,2){$p,S^{(2)}p$}
\end{overpic}\hfill
\begin{overpic}[width=0.48\textwidth]{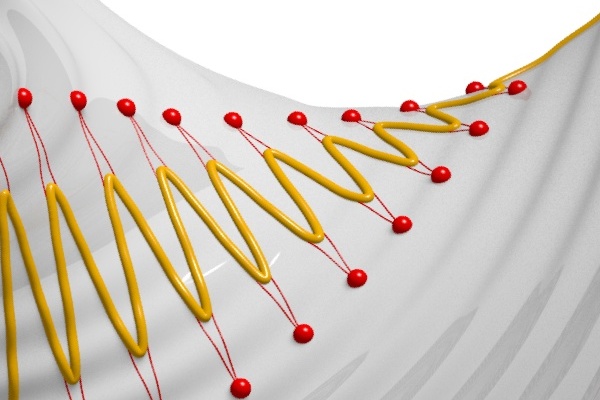}
	\lput(98,2){$p,(S^{(2)})^5 p$}
\end{overpic}
\caption{Geodesic corner-cutting rules are among those where convergence
is not difficult to show. These images show Chaikin's rule
$S^{(2)}$, with the original data in red, and the result of
subdivision as a  yellow geodesic polygon.}
\end{figure}

$S$ is a binary ``dual'' subdivision rule
with combinatorial symmetry about edges. Thus it makes sense to 
define a Riemannian version of $S$ by means of 
averages w.r.t.\ geodesic midpoints of $p_i,p_{i+1}$ as base points, cf.\ Equ.\ \eqref{jw:eq:basepoint:manif}. Using 
$m_{p_i,p_{i+1}}= p_i\oplus {1\over 2}(p_{i+1}\ominus p_i)$, we let
	\begin{align*}
	\Tp_{2i} &= p_i,
	&\Tp_{2i+1}&= m_{p_i,p_{i+1}} \oplus 
		\Big(-\omega(p_{i-1}\ominus m_{p_i,p_{i+1}})
			-\omega(p_{i+2}\ominus m_{p_i,p_{i+1}})
			\Big).
	\end{align*}
 The distance of successive points $\Tp_{2i}$ and $\Tp_{2i+1}$ is bounded
by half the geodesic distance of $p_i,p_{i+1}$ plus
the length of the vector added to the midpoint
in the previous formula. This yields 
the inequality $\delta(Tp)\le {1\over 2}\delta(p) + 2|\omega| {3\over 2}\delta(p)
= ({1\over 2}+3|\omega|)\delta(p)$. Lemma~\ref{jw:lemma1} thus
shows convergence, if $|\omega|<1/6$. 

We cannot easily extend this 
``manifold'' four-point rule to more
general metric spaces. The reason is that we used the linear
structure of the tangent space. A general
discussion of univariate interpolatory rules is found in
\cite{jw:wallner-2014-int}.
\end{Example}

\subsection{Probabilistic interpretation of subdivision
	in metric spaces}

O.~Ebner in \cite{jw:ebner-2011-is,jw:ebner-2014-sars} gave
a probabilistic interpretation of subdivision. This goes as follows.
Consider a linear
subdivision rule as in \eqref{jw:eq:linearrule}, namely
	\begin{align}
	\Sp_i = \sum\nolimits_{j} a_{i-2j} p_j =
	 \avg(a_{i-2j};p_j),
	\quad
	\text{where}\
	a_i\ge 0,\ \sum\nolimits_j a_{i-2j}=1,
	\label{jw:sec:rule:binary}
	\end{align}
acting on data $p\colon \ZZ^s\to\RR^d$.
Consider a stochastic process
$J_0,J_1,\ldots$ defined as the random
walk on $\ZZ^s$ with transition probabilities
	\begin{align*}
	\PP(J_{n+1}\mathord{=}j \mid J_n\mathord{=}i)=a_{i-2j}.
	\end{align*}
 Then the expected value of
$p_{J_{n+1}}$, conditioned on $J_n=j$ is given by
	\begin{align}
	\EE(p_{J_{n+1}} \mid J_n\mathord{=}j) = \Spj,
	\label{jw:eq:expectation}
	\index{expectation, conditional}
	\end{align}
by definition of the expected value. Now the expectation 
$\EE(X)$ of an $\RR^d$-valued random variable $X$
has a characterization via distances: $\EE(X)$
is that constant $c\in\RR^d$ which is closest to $X$ in the sense of 
$\EE(d(X,c)^2)\to\min$. A similar characterization works for the
conditional expectation $\EE(X|Y)$ which is the random variable
$f(Y)$ closest to $X$ in the $L^2$ sense.
These facts inspired a theory of random variables with values in
Hadamard metric spaces developed by
K.-T.\ Sturm \cite{jw:sturm02,jw:sturm03}. The minimizers mentioned
above can be shown to still exist 
if $\RR^d$ is replaced by $\MM$.

Since the way we compute subdivision by Fr\'echet means
is compatible with the distance-based formula for expected values,
Equation \eqref{jw:eq:expectation} holds true also in the case that
both the expectation and the subdivision rule are interpreted
in the Hadamard space sense.
On that basis, O.\ Ebner could
show a remarkable statement on convergence of subdivision rules:

	\index{Hadamard space}
\begin{Satz}{Theorem} {\rm \cite[Th.\ 1]{jw:ebner-2014-sars}}
Consider a binary subdivision rule $\Tp_i=\avg_F(a_{i-2j};p_j)$
with nonnegative coefficients $a_i$. It produces continuous limits
for any data $p_j$ in any Hadamard space $\MM$ if and only if it produces
a continuous
limit function when acting on real-valued data.
\end{Satz}

\begin{proof}[Sketch of proof] With the random walk $(J_i)_{i=0,1,\ldots}$
defined above, \eqref{jw:eq:expectation} directly implies 
	\begin{align}
	(T^n p)_{J_0} \nonumber
	&= \EE(T^{n-1} p_{J_1}\mid J_0)
	= \EE(\EE(T^{n-2} p_{J_2}\mid J_1)\mid J_0) = \ldots
	\\&= \EE(\ldots \EE(\EE(p_{J_n}\mid J_{n-1})\mid J_{n-2})\ldots \mid J_0).
	\label{jw:iterated:conditioning}
	\end{align}
Unlike for $\RR^d$-valued random variables,
there is no tower property for iterated conditioning,
so in general $(T^n p)_{J_0} \ne
\EE(p_{J_n}|J_0)$. That expression has a different interpretation:
$T$ is analogous to the linear rule $S$ of 
\eqref{jw:sec:rule:binary}, which is nothing but the restriction of
the general rule $T$ to data in Euclidean spaces. Its $n$-th power
$S^n$ is a linear rule of
the form $(S^n q)_i = \sum a^{[n]}_{i-2^n j}q_j$, and we have
	\begin{align}
	\EE(q_{J_n}\mid J_0) = (S^n q)_{J_0},
	\quad \text{if $S$ acts linearly on data $q\colon \ZZ^s\to\RR^d$.}
	\label{jw:iterated:linear}
	\end{align}
 This follows either directly (computing the coefficients
of the $n$-th iterate $S^n$ corresponds to computing
transition probabilites for the $n$-iterate of the random walk),
or by an appeal to the tower property in 
\eqref{jw:iterated:conditioning}.

Sturm \cite{jw:sturm02} showed
a Jensen's inequality for continuous convex functions
$\Psi$, 
	\index{Jensen's inequality}
	\begin{align*}
	\Psi\Big(\EE(\ldots (\EE(p_{J_n}\mid J_{n-1})\ldots \mid J_0)\Big)
	\le \EE\Big(\Psi(p_{J_n})\mid J_0\Big).
	\end{align*}
We choose $\Psi=\dMM(\cdot ,x)$ and observe that $q_{J_n}=\dMM(p_{J_n},x)$ is a
real-valued random variable. Combining Jensen's inequality with
\eqref{jw:iterated:conditioning} and \eqref{jw:iterated:linear} yields
	\begin{align*}
	\dMM(T^n p_i,x)& \le 
		\sum\nolimits_k a^{[n]}_{i-2^nk} \dMM(p_k,x),
		&& (\text{for any}\ x)
	\\
	\dMM(T^n p_i,T^n p_j) 
	& \le \sum\nolimits_{k,l} 
		a^{[n]}_{i-2^nk}
		 a^{[n]}_{j-2^nl}
		 \dMM(p_k,p_l) 
	&& (\text{by recursion)}.
	\end{align*}
To continue, we need some information 
on the coefficients $a^{[n]}_i$. For that, we use the 
limit function $\phi\colon\RR^s\to[0,1]$ generated by applying
$T$ (or rather, $S$), to the
delta sequence. By construction (see Lemma \ref{jw:lemma1}),
$|a^{[n]}_j-\phi(2^{-n} j)|\to 0$ as $n\to\infty$.
These ingredients allow us to show 
existence of $n$ with $T^n$ contractive.
\end{proof}

As a corollary we get, for instance, that subdivision with nonnegative
coefficients works in $\Pos_n$ in the same way as in linear spaces,
as far as convergence is concerned. Since $\Pos_n$ is not only a
Hadamard metric space, but even a smooth Riemannian manifold, also
the next section will yield a corollary regarding $\Pos_n$.

\subsection{The convergence problem in manifolds}
\label{jw:ss:subdivmanif}

The problem of convergence of subdivision rules in manifolds (Riemannian
manifolds, groups, and symmetric spaces) was at first treated by means
of so-called proximity inequalities which compare linear rules with
their analogous counterparts in manifolds. This approach was successful
in studying smoothness of limits (see Section \ref{jw:sec:smoothness}
below), but less so for convergence. Unless subdivision rules are of a 
special kind (interpolatory, corner-cutting$\,,\ldots$)
convergence can typically be shown
only for ``dense enough'' input data, with very small
bounds on the maximum allowed density. On the other hand numerical
experiments demonstrate that a manifold rule analogous to a convergent
linear rule usually converges. This
discrepancy between theory and practice is of course unsatisfactory 
from the viewpoint of theory, but is not so problematic from the viewpoint
of practice. The reason is the stationary nature of subdivision --- if
$\delta(p)$ is too big to infer existence of a continuous limit
$S^\infty p$, we can check if $\delta(\Skp)$ is small enough instead.
As long as $S$ converges, this leads to an a-posteriori proof
of convergence.

More recently, convergence of subdivision rules of the form
$\Sp_i=\avg_F(a_{i-Nj};p_j)$ in Riemannian manifolds has been 
investigated along the lines of Lemma~\ref{jw:lemma1}.
This work is mainly based on the methods of
H.~Karcher's seminal paper \cite{jw:karcher}. So far, only the
univariate case of data $p\colon\ZZ\to\MM$ has been treated
successfully, cf.\
\cite{jw:wallner-2011-sym,jw:huening-subdiv-2017,jw:huening-2019}.

There are two main cases to consider. In Cartan-Hadamard manifolds
(curvature $\le 0$) the Fr\'echet mean is well defined
and unique also if weights are allowed to be negative
\cite[Th. 6]{jw:huening-subdiv-2017}. Subdivision
rules are therefore globally and uniquely defined. 
We have the following result:

\begin{Satz}{Prop.} {\rm \cite[Th.\ 11]{jw:huening-subdiv-2017}}
Consider a univariate subdivision rule
$\Sp_i=\avg_F(a_{i-Nj};p_j)$ acting on sequences in 
a Cartan-Hadamard manifold $\MM$. Consider also
the norm $\|S^*\|$ of its linear derived subdivision rule
according to Example~\ref{jw:ex:linear}. If
	\begin{align*}
	\gamma = {1\over N}\|S^*\|
	< 1,
	\end{align*} 
 then $S$ meets the conditions of Lemma~\ref{jw:lemma1} (with
contractivity factor $\gamma$) and 
produces continuous limits. 
\end{Satz}

This result is satisfying because it allows us to infer convergence
from a condition which is well known in the linear case, cf.\
\cite{jw:dyn:1992:SSIC}. If ${1\over N}\|S^*\|\ge 1 $, we can instead
check if one of
${1\over N^n}\|S^{*n}\|$, $n=2,3,\ldots$
is smaller than $1$. If this is the case,
then the manifold subdivision rule analogous to the linear rule $S^n$
converges.

\begin{figure}[t]
\centering\unitlength0.0098\textwidth
\begin{picture}(60,33)\put(-6,0){
\put(0,0){\includegraphics[width=30\unitlength]{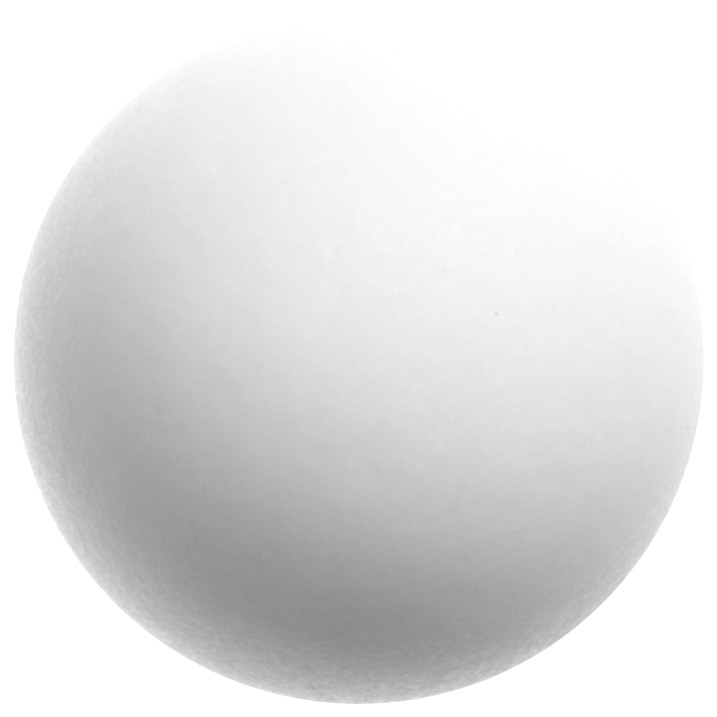}}
\lput(60,0){\includegraphics[width=30\unitlength]{sphere1.jpg}}
\put(0,0){\includegraphics[width=30\unitlength]{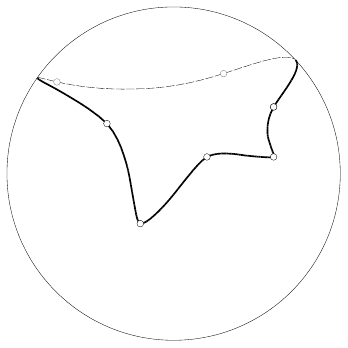}}
\lput(60,0){\includegraphics[width=30\unitlength]{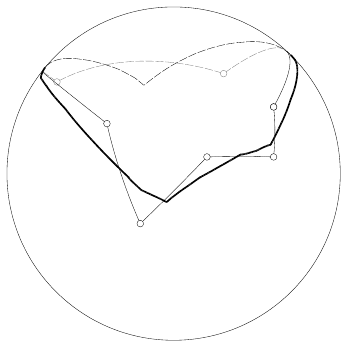}}
\lput(5,1){\scriptsize${1\over 16}(-1,0,9,1,9,0,-1)$}
\put(56,1){\scriptsize${1\over 32}(-1,-1,21,13,13,21,-1,-1)$}
\cput(15,4){$\delta(p)<0.31$}
\cput(45,4){$\delta(p)<0.4$}
}\end{picture}

\caption{Subdivision rules $\Spj=\avg_F(a_{j-2i};p_i)$
based on the Fr\'echet mean operating on
sequences on the unit sphere. The images visualize the interpolatory
4-point rule (left) and a rule without any special properties.
We show the coefficient sequence $a_j$ and the bound on $\delta(p)$
which ensures convergence.}
\label{jw:fig:spheres}
\end{figure}

\Paragraph{Subdivision in Riemannian manifolds with positive curvature}
Recent work \cite{jw:huening-2019} deals with spaces of positive
curvature, and initial results have been achieved on the unit sphere, for
subdivision rules of the form $\Sp_i=\avg_F(a_{i-2j};p_j)$. 
Figure~\ref{jw:fig:spheres} shows two examples. One aims
at finding a bound $\delta_0$ such that for all
data $p$ with $\delta(p)<\delta_0$, $S$ acts in a contractive way so that
Lemma~\ref{jw:lemma1} shows convergence.

Rules defined in a different way are sometimes much easier to analyze.
E.g.\ the Lane-Riesenfeld subdivision rules defined by midpoint insertion,
followed by $k$ rounds of averaging, can be transferred to any complete
Riemannian manifold as a corner-cutting rule and will enjoy continuous
limits, see Example~\ref{jw:ex:cornercutting}. Similarly, the
interpolatory four-point rule can be generalized to the manifold
case in the manner described by Example~\ref{jw:ex:interpolatory}, and
will enjoy continuous limits. The generalization via the Fr\'echet mean
(Fig.~\ref{jw:fig:spheres}) on the other hand, is not so easy to analyze.
The approach by \cite{jw:huening-2019} is to control $\delta(\Sp)$ 
by introducing a family $S^{(t)}$, $0\le t\le 1$, of rules where $S^{(0)}$ is 
easy to analyze, and $S^{(1)}=S$. If one manages to show 
$\delta(S^{(0)}p)<\gamma_1 \delta(p)$ and
$\|{d\over dt}S^{(t)}p_i\|\le C\delta(p)$, then the length
of each curve $t\mapsto S^{(t)} p_i$ is bounded by $C\delta(p)$, and
	\begin{align*}
	\delta(\Sp) 
	&\le
	\sup_i\dMM(\Sp_{i},S^{(0)}p_{i}) 
	+ \delta(S^{(0)}p)
	+\sup_i\dMM(\Sp_{i+1},S^{(0)}p_{i+1}) 
	\\& \le (\gamma_1+2C)\delta(p).
	\end{align*}
Contractivity is established if $\gamma_1+2C<1$, in which case 
Lemma~\ref{jw:lemma1} shows convergence.
The bounds mentioned in Fig.~\ref{jw:fig:spheres} have been found in 
this way. Estimating the norm of the derivative mentioned above involves 
estimating the eigenvalues of the Hessian of 
the right hand side of \eqref{jw:eq:frechetmean}.

\Paragraph{The state of the art regarding convergence of
refinement schemes}
Summing up, convergence of geometric subdivision rules
is treated in a satisfactory manner 
for special rules (interpolatory, corner-cutting), for rules in special
spaces (Hadamard spaces and Cartan-Hadamard manifolds), and
in the very special case of the unit sphere and univariate rules. 
General manifolds with positive curvature have not been treated.
Multivariate data are treated only in Hadamard metric spaces and
for subdivision rules with nonnegative coefficients. In other situations,
we know that convergence happens only for ``dense enough'' input data,
where the required theoretical
upper bounds on $\delta(p)$ are very small compared
to those inferred from numerical evidence.

\section{Smoothness analysis of subdivision rules}
\label{jw:sec:smoothness}

For linear subdivision rules, the question of smoothness of limits
can be considered as largely solved, 
the derived rule $S^*$ introduced in Example~\ref{jw:ex:linear} being the
key to the question if limits are smooth.
Manifold subdivision rules do not always enjoy the same smoothness as the
linear rules they are derived from.
The constructions
mentioned in Section~\ref{jw:sec:subdivision} basically yield manifold 
rules whose limits enjoy $C^1$ resp.\ $C^2$ smoothness if the original
linear rule has this property, but this general statement is no longer true
if $C^3$ or higher smoothness is involved. Manifold
rules generated via Fr\'echet means or via projection 
\cite{jw:grohs-2009-ps,jw:xieyu-smoothnessequiv} retain the smoothness
of their linear counterparts. Others, e.g.\ 
constructed by means of averages w.r.t.\
basepoints in general do not. This is 
to be expected, since the choice of basepoint introduces an element
of arbitrariness into manifold subdivision rules. 
The following paragraphs discuss the method of {\em proximity
inequalities} which was successfully employed in treating the smoothness
of limits.

\subsection{Derivatives of limits}

A subdivision rule $S$ acting on a sequence $p$ in $\RR^d$ converges
to the limit function $S^\infty p$, if the refined data $\Skp$,
interpreted as samples of functions $f_k$ at the finer grid $N^{-k}\ZZ$,
approach that limit function (see Definition~\ref{jw:defn:convergence}):
	\begin{align*}
	(S^\infty p)(\xi)
	 & \approx 
		f_k(\xi)=(\Skp)_{N^k\xi}, 
	\end{align*}
whenever $N^k\xi$ is an integer. 
A similar statement holds for derivatives, which are approximated
by finite differences. With $h=N^{-k}$, we get
	\begin{align*}
	(S^\infty p)'(\xi)
	& \approx 
	{f_k(\xi+h)-f_k(\xi)\over h}
	= N^k((\Skp)_{N^k\xi+1}-(\Skp)_{N^k\xi})
	\\&= (\Delta (NS)^k p)_{N^k\xi}
	=
	 ( S^{*k}\Delta p) _{N^k\xi}
	.
	\end{align*}
Here $S^*$ is the derived rule defined by the
relation $S^*\Delta=N\Delta S$,  see Ex.~\ref{jw:ex:linear}.
For the $r$-th derivative of the limit function we get
	\begin{align*}
	(S^\infty p)^{(r)}(\xi)
	& \approx 
	 (\Delta^r (N^rS)^k p)_{N^k\xi}
	=
	\def\TMP{{\unitlength1ex\begin{picture}(0.01,0.1)\tiny
		\cput(0.2,1.6){$r$ times}
		\cput(0.2,0.4){$\overbrace{\hphantom{\scriptsize **\cdots*}}$}
		\end{picture}}}
	((S^{**\TMP\cdots *})^k\Delta^r p)_{N^k\xi}
	.
	\end{align*}
These relations, except for references to derived rules,
are valid even if $S$ does not act linearly.
$S$ could be a manifold rule expressed in 
a coordinate chart, or it could be acting on a surface contained in $\RR^d$.

If $S$ does act linearly, 
one proves that $S$ has $C^1$ smooth limits, 
if $S^*$ has continuous ones, and in that case
$(S^\infty p)' = S^{*\infty}\Delta p$. To treat
higher order derivatives, this statement can be iterated. 
For multivariate data $p_i$, $i\in\ZZ^s$, the situation is analogous
but more complicated to write down. 
For the exact statements, see
\cite{jw:dyn:1992:SSIC,jw:cavaretta-1991-ss}.

\subsection{Proximity inequalities}
\label{jw:sec:proximity}
\index{proximity inequality}

\Paragraph{Smoothness from proximity}
Manifold subdivision rules were first systematically analyzed with regard to
derivatives by \cite{jw:wallner-2004-cca}.
The setup is a linear rule $S$ and a nonlinear rule $T$ both acting
on data contained in the same space $\RR^d$.
$T$ could be a manifold version of $S$, with $\RR^d$ being a coordinate chart
of the manifold; or $T$ could act on points of a surface contained in
$\RR^d$. Then $S$, $T$ are in proximity, if 
	\begin{align}
	\sup\nolimits_i \|\Sp_i-\Tp_i\| \le C\delta(p)^2.
	\label{jw:eq:proximity0}
	\end{align}
This formula is motiviated by a comparison of the shortest
path between two points within in a surface (which is a geodesic segment),
with the shortest path
in Euclidean space (which is a straight line). These two paths differ
by exactly the amount stated in \eqref{jw:eq:proximity0}.
Two statements 
were shown in \cite{jw:wallner-2004-cca}:

\begin{list}{}{\itemsep-\parsep}

\item[(1)] Certain manifold subdivision rules 
$T$ derived from a convergent linear rule $S$ obey 
the proximity inequality \eqref{jw:eq:proximity0}
whenever data are dense enough (i.e., $\delta(p)$ is small enough).

\item[(2)] in that case, if limit curves of $S$ 
enjoy $C^1$ smoothness, then $T$ produces continuous limit curves for data
with $d(p)$ small enough; and all continuous limit curves enjoy $C^1$
smoothness.

\end{list}

\noindent
To demonstrate how proximity inequalities work, we prove a convergence
statement like the ones given by \cite[Th.\ 1]{jw:wallner-2006-is} or
\cite[Th.\ 2+3]{jw:wallner-2004-cca} (with slightly different proofs).

\begin{Satz}{Prop.}
Assume the setting of Equ.\ \eqref{jw:eq:proximity0}, with a subdivision
rule $T$ being in proximity with a linear subdivision rule $S$.
We also assume ${1\over N}\|S^*\|<1$.\footnote{implying convergence
of the linear rule $S$. $N$ is the dilation factor, $S^*$ is the derived
rule, cf.\ Ex.~\ref{jw:ex:linear}.}
Then $T$ produces continuous limit curves from data
$p$ with $\delta(p)$ small enough.
\end{Satz}

\begin{proof}
Generally $ \sup_i\|p_i-q_i\| \le K \implies 
\delta(p)\le \delta(q)+2K$. Thus 
\eqref{jw:eq:proximity0} implies
		\begin{align*}
		\delta(\Tp) 
		\le \delta(\Sp)+2C \delta(\Tp)^2
		\le {1\over N}\|S^*\| \delta(p)+2C \delta(p)^2.
		\end{align*}
Choose $\epsilon>0$ with 
$\lambda :={1\over N}\|S^*\|+2C\epsilon <1$. If $\delta(p)<\epsilon$,
then $T$ is contractive:
	\begin{align*}
	\delta(\Tp)
	\le ({1\over N}\|S^*\|+2C\delta(p))\delta(p) 
	\le\lambda\delta(p)
	.
	\end{align*}
By recursion, $\delta(T^{k+1}p)\le\lambda\delta(T^k p)$.
As to the displacement-safe condition of Lemma~\ref{jw:lemma1},
recall from Example~\ref{jw:ex:linear} that
$S$ has it. For $T$, observe that 
	\begin{align*}
	\|\Tp_{Ni}-p_i\|
	\le\|\Tp_{Ni}-\Sp_{Ni}\|+\|\Sp_{Ni}-p_i\|
	\le C\delta(p)^2+C'\delta(p)
	\le (\epsilon C+C')\delta(p).
	\end{align*}
Now Lemma~\ref{jw:lemma1} shows convergence.
\end{proof}

\noindent
The convergence of vectors $N^k\Delta T^k p$ to 
derivatives of the limit function $T^\infty p$
is proved in a way which is analogous in principle. The method was extended
to treat $C^2$ smoothness 
by \cite{jw:wallner-2005-sass}, 
using the proximity condition 
	\begin{align*}
	\sup\nolimits_i \|\Delta\Sp_i-\Delta\Tp_i\| 
		\le C(\delta(p)\delta(\Delta p)
		+\delta(p)^3)
	.
	\end{align*}
A series of publications treated $C^2$ smoothness
of Lie group subdivision rules based on log/exp averages 
\cite{jw:wallner-2007-sl,jw:grohs-2008-lgsd}, 
the same in symmetric spaces
\cite{jw:wallner-2011-sym}, $C^1$ smoothness
in the multivariate case \cite{jw:grohs-2008-sbg}, higher order
smoothness of interpolatory rules in groups
\cite{jw:grohs-2008-mlie,jw:xieyu-lie}, and higher order smootheness
of projection-based
rules \cite{jw:grohs-2009-ps,jw:xieyu-smoothnessequiv}. The proximity
conditions involving higher order smoothness become rather complex, 
especially in the multivariate case.

\Paragraph{Smoothness equivalence}
If a manifold subdivision rule $T$ is created on basis of a linear rule $S$,
it it interesting to know if
the limit functions of $T$ enjoy the same smoothness as the limits of $S$.
For $C^1$ and $C^2$ smoothness, $T$ can basically
be constructed by any of the methods described above, and it will enjoy
the same smothness properties as $S$ (always assuming that convergence happens,
and that the manifold under consideration is itself as smooth as the
intended smoothness of limits).
This smoothness equivalence breaks down for $C^k$ with $k\ge 3$. 

A manifold subdivision rule based on the 
log/exp construction, using averages w.r.t.\ basepoints, 
		\begin{align*}
		\Tp_i = \avg_{m_i}(a_{i-Nj};p_j), 
		\index{log/exp subdivision}
		\end{align*}
does not enjoy $C^k$ smoothness equivalence for $k\ge 3$
unless the base points $m_i$
obey a technical condition which can be satisfied e.g.\ if
they themselves are produced by certain kinds of subdivision
\cite{jw:yu-cksmoothness,jw:grohs-2010-gpa}.  Necessary
and sufficient conditions for smoothness equivalence are
discussed by \cite{jw:duchamp-smothnessequiv-2015}.
We pick one result whose proof is based on this method (using \eqref{jw:eq:self}
for a ``base point'' interpretation of Fr\'echet means):

\index{Frechet mean@Fr\'echet mean}
\begin{Satz}{Theorem}\label{jw:th:philipp} 
{\rm \cite[Th.\ 4.3]{jw:grohs-2010-gpa}}
Let $S$ be a stable\footnote{``Stable'' means existence of constants 
$C_1,C_2$ with $C_1\|p\|\le \|S^\infty p\|_\infty
\le C_2\|p\|$ for all input data where $\|p\|:=\sup_i\|p_i\|$ is bounded.
Stable rules with $C^n$ limits generate
polynomials of degree $\le n$, which is a property used in the proof.}
subdivision rule $\Sp_i=\avg(a_{i-Nj}$; $p_j)$
acting on data $p\colon\ZZ^s\to\RR^d$,
which is convergent with $C^n$ limits. 
Then all continuous limits of its
Riemannian version $\Tp_i$ $=$ $\avg_F(a_{i-Nj}$; $p_j)$ likewise
are $C^n$.
\end{Satz}

We should also mention that proximity conditions relevant to the
smoothness analysis of manifold subdivision rules can take 
various forms, cf.\ the ``differential'' proximity condition of
\cite{jw:grohs-2010-st,jw:duchamp-smothnessequiv-2015,jw:yu:newproximity}.

Finally we point out a property which manifold rules share with
linear ones: For any univariate linear rule $S$ which has $C^k$ limits,
the rule $A_{1/2}^k\circ S$ has limits of smoothess $C^{n+k}$,
where $A_{1/2}$ is midpoint-averaging as described by
Ex.~\ref{jw:ex:cornercutting}.
It has been shown in \cite{jw:yu-smoothing-2018}
that an analogous statement
holds true also in the manifold case, for a general class of
averaging operators.
	\index{corner cutting}

\subsection{Subdivision of Hermite data}
	\index{Hermite subdivision}
	\label{jw:sec:hermite}

\def\PV{\mathord{\textstyle{p\choose v}}}
\def\QV{\mathord{\textstyle{q\choose w}}}

\begin{figure}[t]
\quad\begin{overpic}[height=.25\textwidth]{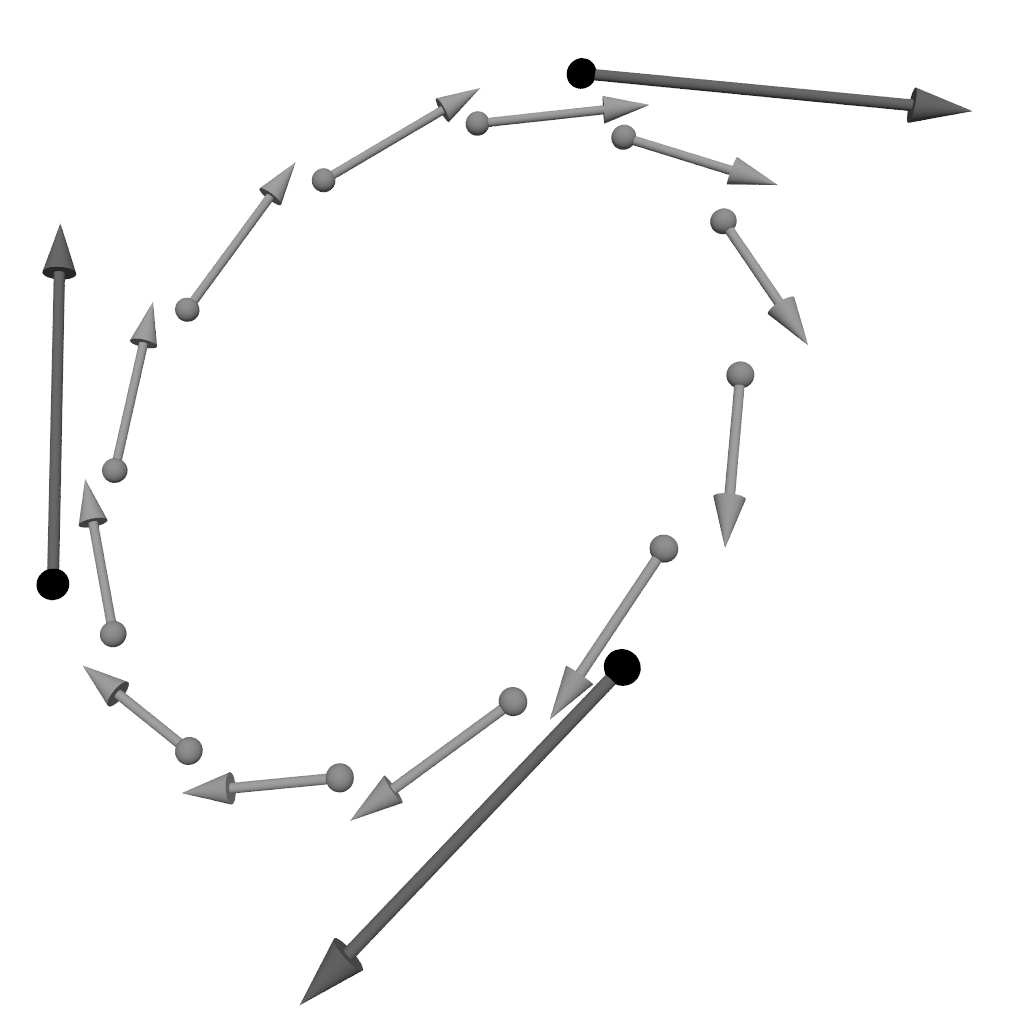}
	\put(50,10){$(p_0,v_0)$}
	\cput(5,80){$(p_1,v_1)$}
	\cput(75,95){$(p_2,v_2)$}
	\end{overpic}
\begin{overpic}[height=.25\textwidth]{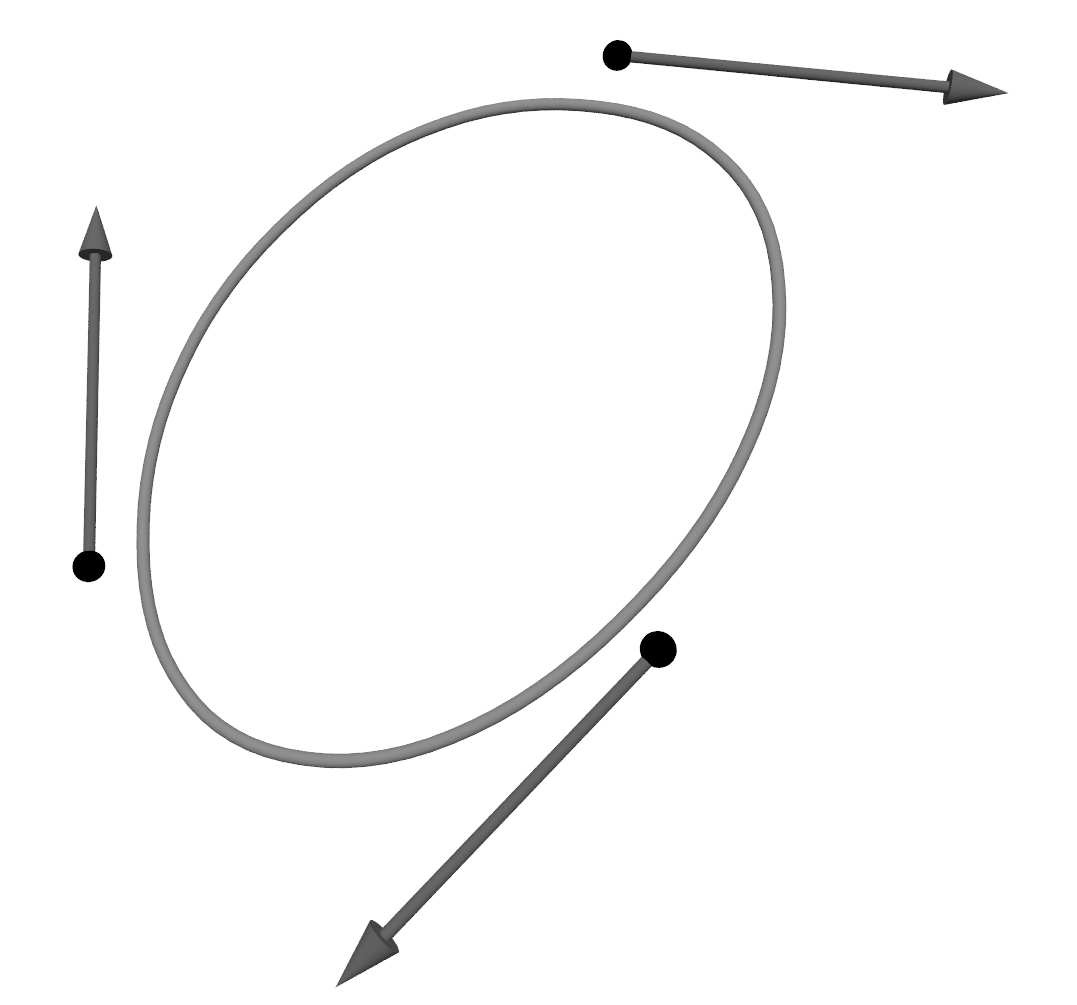}
	\cput(30,30){$f$}
	\end{overpic}\hfill
\includegraphics[height=.3\textwidth]{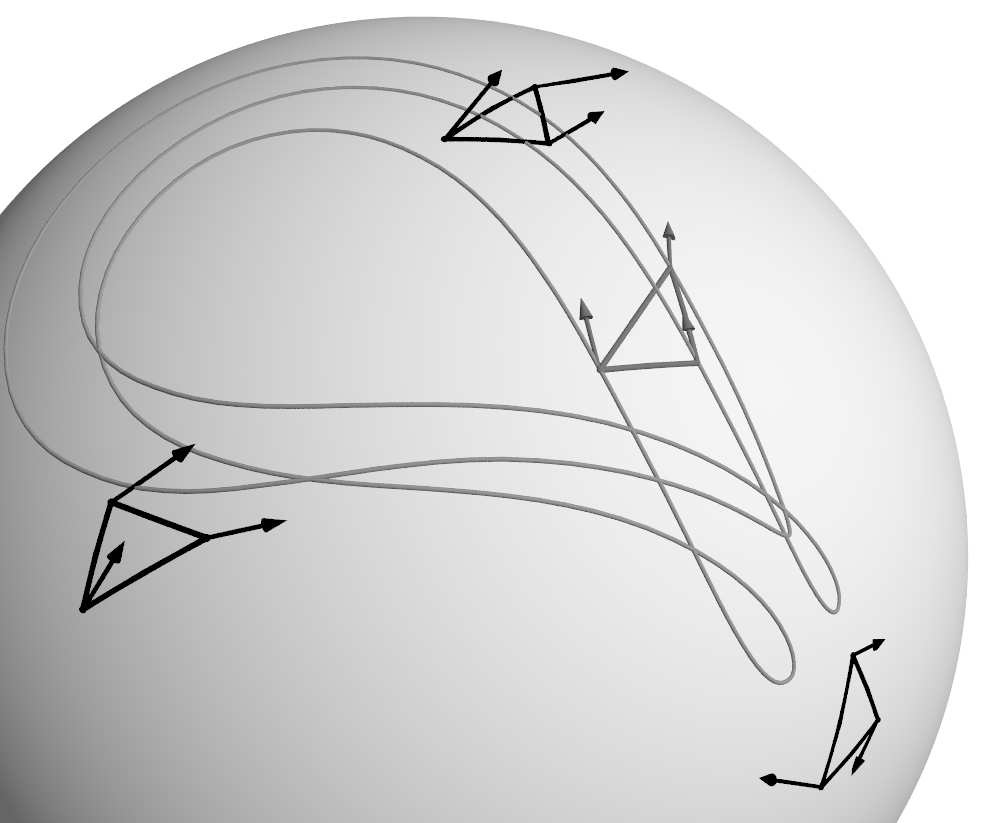}
\caption{{\it Left:} Hermite data $(p_i, v_i)$ in $\RR^2$ and
the result of one round subdivision by a linear Hermite rule $S$.
{\it Center:} Limit curve $f$ ($f'$ is not shown).
{\it Right:} Hermite data $(p_i, v_i)$ in the group $\SO_3$,
and the limit curve generated by a group version of $S$. Points
$p_i\in\SO_3$ and tangent vectors $v_i\in T_{p_i}\SO_3$ are visualized
by means of their action on a spherical triangle. These figures
appeared in \protect\cite{jw:moosmueller-hermite2}
(reprinted with permission).}
\end{figure}

Hermite subdivision is a refinement process acting not on points, but
on tangent vectors, converging to a limit and its derivative 
simultaneously.
In the linear case, data $(p,v)\colon\ZZ\to\RR^d
\times\RR^d$ undergo subdivision by a rule $S$ which obeys basic
shift invariance $SL=L^NS$. The interpretation of
$p_i$ as points and $v_i$ as vectors leads to 
	\begin{align}
	S\PV_i = 
	\bigg(\begin{array}{l}
	 \sum_j a_{i-Nj} p_j + \sum_j b_{i-Nj} v_j \\
	 \sum_j c_{i-Nj} p_j + \sum_j b_{i-Nj} v_j. 
	\end{array}\bigg)\,, \quad
	\text{where}\ 
	\begin{cases}
	 \sum_j a_{i-Nj}=1,
	\\ \sum_j c_{i-Nj}=0.
	\end{cases}
	\label{jw:eq:hermite}
	\end{align}
$S$ is invariant w.r.t.\ 
translations, which act via $p\mapsto p+x$ on points,
but act identically on vectors. Iterated refinement creates data
$S^k\PV$ converging to a limit 
$f\colon\RR\to\RR^d$,
	\begin{align*}
	\bigg({f(\xi)\atop f'(\xi)}\bigg) = 
	\lim_{k\to\infty} 
	\bigg({1\atop 0}\ {0 \atop N^k}\bigg)
	S^k \bigg({p\atop v}\bigg)_{N^k\xi}\ , \quad
	\text{whenever}\ N^k\xi\in\ZZ.
	\end{align*}
We say that $S$ converges, if the limit $(f,f')$ exists and $f$ enjoys
$C^1$ smoothness, with $f'$ then being continuous.
A manifold version of $S$, operating on data
	\begin{align*}
	\PV\colon \ZZ\to T\MM, \quad \text{i.e.,}\ v_i \in T_{p_i}\MM,
	\end{align*}
	\index{parallel transport}
faces the difficulty
that each $v_i$ is contained in a different vector space. One
possibility to overcome this problem is to employ parallel transport 
$\PT_p^q\colon T_p\MM\to T_q\MM$ between tangent spaces. 
In Riemannian manifolds, a natural choice for $\PT_p^q$ 
is parallel transport w.r.t.\ the canonical Levi-Civita connection
along the shortest geodesic
connecting $p$ and $q$, cf.\ \cite{jw:docarmo}.
In groups, we can simply choose 
$\PT_p^q$ as left translation by $qp^{-1}$ resp.\ the differential
of this left translation. Then the definition
	\begin{align*}
	S\PV = \QV \ \text{with}\
	\bigg\{\begin{array}{*9{c@{\,}}}
	q_i & = &
		m_i\oplus
		& \big(\sum_j a_{i-Nj} (p_{j}\ominus m_i)
		&+ \sum_j b_{i-Nj} 
		\PT_{p_{j}}^{m_i} v_{j}\big)
	\\[0.5ex]
	w_i & = &
		\PT_{m_i}^{q_i}
		& \big(\sum_j c_{i-Nj} (p_{j}\ominus m_i)
		&+ \sum_j d_{i-Nj} 
		\PT_{p_{j}}^{m_i} v_{j}\big)
	\end{array}
	\end{align*}
is meaningful (provided 
the base point $m_i$ is chosen close to $p_{\lfloor i/N\rfloor}$).
In a linear space, this expression reduces to \eqref{jw:eq:hermite}.
C.\ Moosm\"uller could show $C^1$ smoothness of limits
of such subdivision rules, by methods 
in the spirit of Section~\ref{jw:sec:proximity},
see \cite{jw:moosmueller-hermite-2016, jw:moosmueller-hermite2}.

\subsection{Subdivision with irregular combinatorics}
\label{jw:sec:irregular}

\begin{figure}[b]
\includegraphics[width=0.245\textwidth]{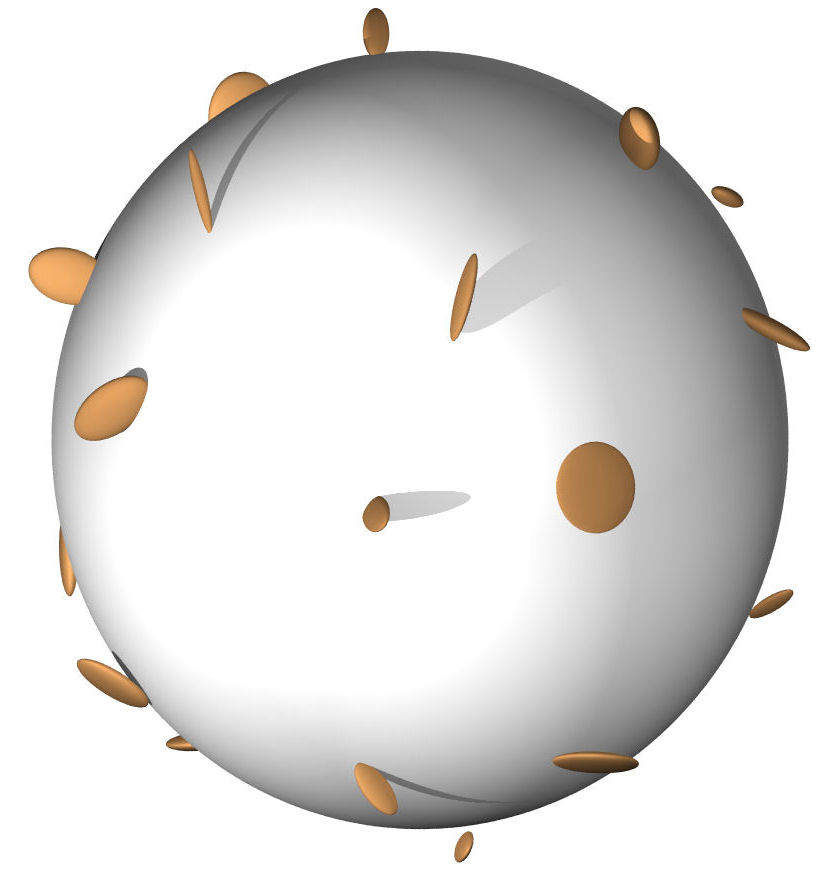}\hfill
\includegraphics[width=0.245\textwidth]{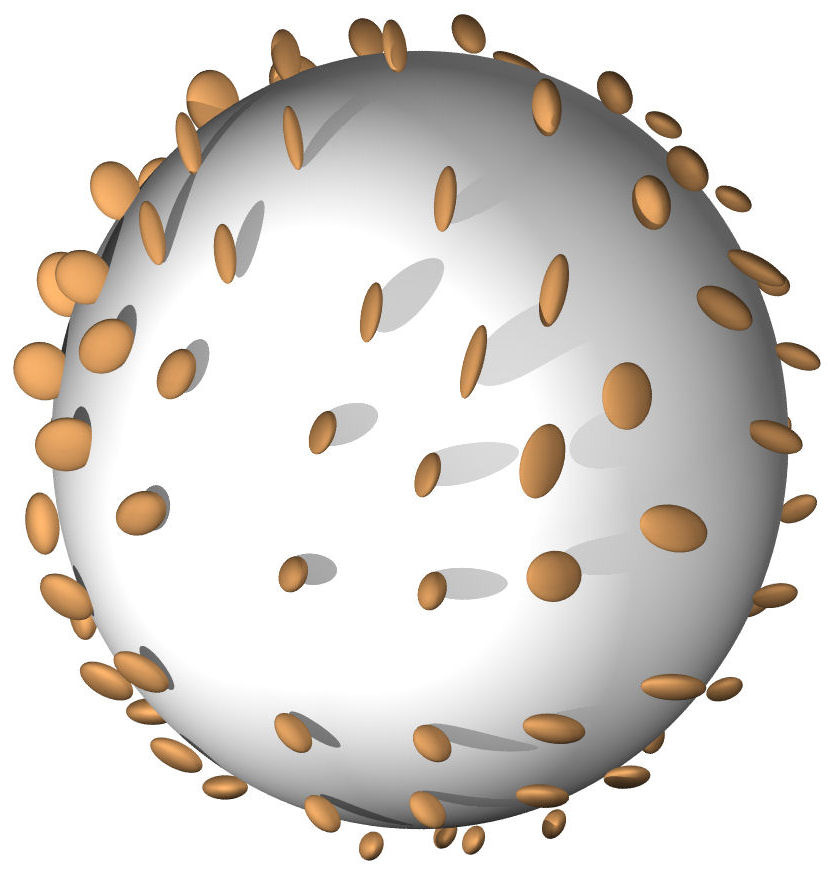}\hfill
\includegraphics[width=0.245\textwidth]{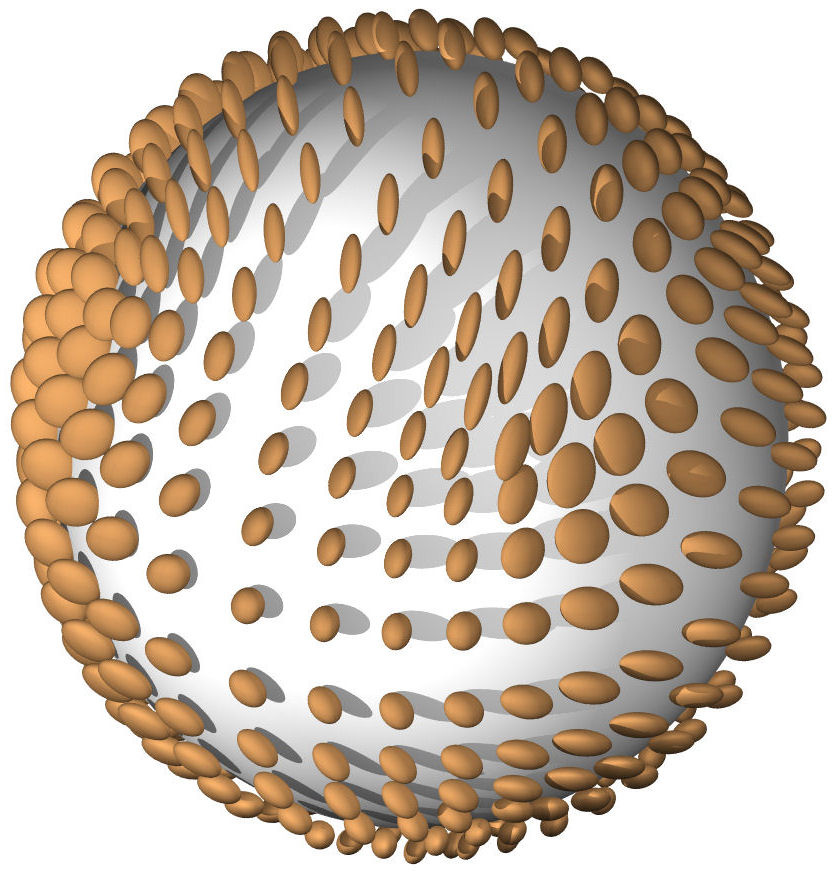}\hfill
\includegraphics[width=0.245\textwidth]{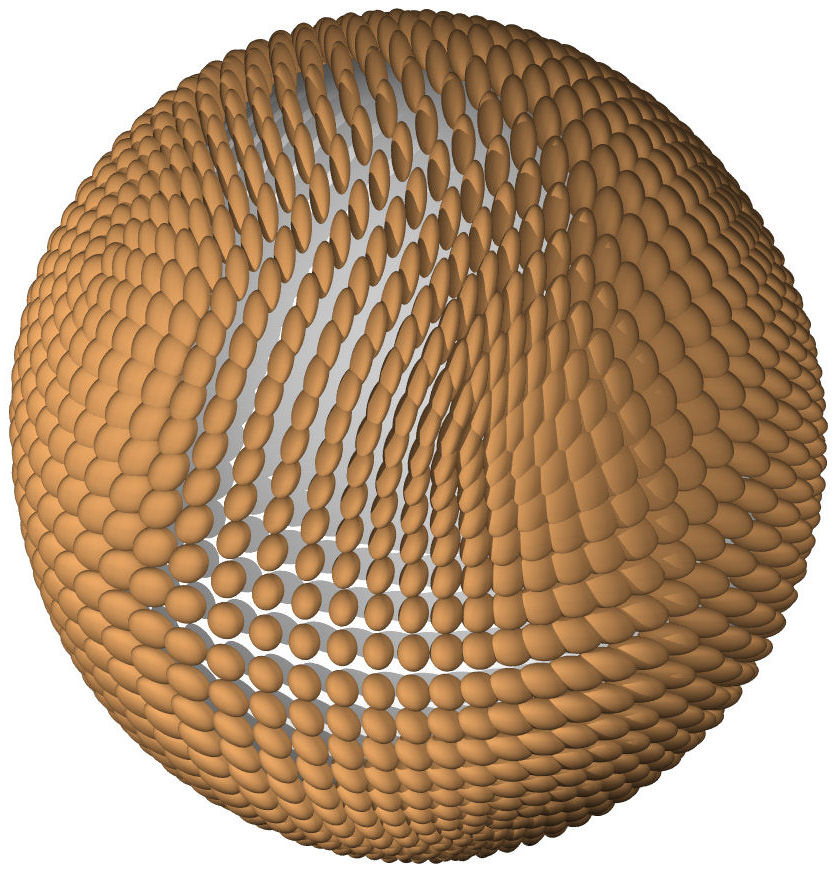}
\caption{Here data $p_i$ in the unit sphere $\Sigma^2$ and
$\Pos_3$-valued data $q_i$ are visualized by placing the ellipsoid with
equation $x^T q_i x=1$ in the point $p_i\in \Sigma^2$. Both data undergo
iterative refinement by means of a Riemannian version $S$ of
the Doo-Sabin subdivision rule. For given initial data $p,q$ which
have the combinatorics of a cube, the four images show
$S^jp$ and $S^j q$, for $q=1,2,3,4$ (from left).
The correspondence $(\Skp)_i\mapsto (S^kq)_i$ converges to a $C^1$
immersion
$f\colon \Sigma^2\to \Pos_3$ as $k\to\infty$.
These figures appeared in \protect\cite{jw:weinmann-2009-nsi}
(reprinted with permission).}
\label{jw:fig:irregular}
\end{figure}

A major application of subdivision is in Computer Graphics,
where it is 
ubiquitously used as a tool to create surfaces from a finite number of
handle points whose arrangement is that of the vertices of
a 2D discrete surface. That surface
usually does not have the combinatorics of a regular grid.

Two well known subdivision rules acting on such
data are the Catmull-Clark rule and the Doo-Sabin rule, see
\cite{jw:catmullclark1978,jw:doosabin1978}. Such
subdivision rules create denser and denser discrete surfaces which
are mostly regular grids but retain a constant number of combinatorial
singularities. This implies that the limit surface is locally
obtained via Definition~\ref{jw:defn:convergence}, but with a nontrivial
overlapping union of several such limits as one approaches
a combinatorial singularity. A systematic way of analyzing convergence
and smoothness was found by U.\ Reif \cite{jw:reif-1995-uasd}, see also the
monograph \cite{jw:Peters2008}. There is a wealth of contributions
to this topic, mostly because of its relevance for Graphics.

A.\ Weinmann in \cite{jw:weinmann-2009-nsi, jw:weinmann-2012-ad, jw:weinmann-2012-mf} studied intrinsic manifold versions of such subdivision rules.
They are not difficult to define, since the linear subdivision rules
which serve as a model are defined in terms of averages. We do not
attempt to describe the methods used for establishing convergence
and smoothness of limits other than to say that a proximity condition
which holds between a linear rule $S$ and a nonlinear rule $T$ eventually 
guarantees that in the limit, smoothness achieved by $S$ carries over
to $T$ --- the perturbation incurred by switching from a linear space
to a manifold is not sufficient to destroy smoothness.
Figure~\ref{jw:fig:irregular} illustrates a result obtained by 
\cite{jw:weinmann-2009-nsi}.

\section{Multiscale transforms}
\index{multiscale transform|(}

\subsection{Definition of intrinsic multiscale transforms}

A natural multiscale representation of data, which does not suffer from 
distortions caused by the choice of more or less arbitrary
coordinate charts, is required to be based on operations which 
are themselves adapted to the geometry of the data. This topic
is intimately connected to subdivision, since upscaling operations
may be interpreted as subdivision.

A high-level introduction of certain kinds of multiscale decompositions
is given by \cite{jw:grohs-2012-ds}. We start with an elementary
example.

\begin{Example}{A geometric Haar decomposition and reconstruction
	procedure}\label{jw:ex:haar}\relax
Consider data
$p\colon\ZZ\to\MM$, and the upscaling rule $S$ and downscaling rule $D$,
	\begin{align*}
	(\ldots,p_0,p_1,\ldots)
		& \stackrel{S}\longmapsto 
		(\ldots,p_0,p_0, p_1,p_1,\ldots)
	\\[-1ex]
	(\ldots,p_0,p_1,\ldots)
		& \stackrel{D}\longmapsto 
		(\ldots,m_{p_0,p_1},m_{p_1,p_2},\ldots), \
	\text{where}\ m_{a,b}=a\oplus{1\over 2}(b\ominus a).
	\end{align*}
The use of $\oplus$ and $\ominus$ refers to the exponential mapping,
as a means of computing differences of points, and adding vectors to points.
$D$ is a left inverse of $S$ but not vice versa: $SDp\ne p$ in general.
However, if 
we store the difference between $p$ and $SDp$ as {\em detail vectors} $q$:
	\begin{align*}
	(\ldots,q_0,q_1,\ldots) =
		(\ldots, p_0\ominus m_{p_0,p_1}, p_2\ominus m_{p_2,p_3},\dots)
	\end{align*}
 then the reconstruction procedure
	\begin{align*}
	p_{2i} &= m_{p_{2i},p_{2i+1}} \oplus q_i, 
	& p_{2i+1} &= m_{p_{2i},p_{2i+1}} \ominus q_i
	\end{align*}
 recovers the information destroyed by downsampling.
\end{Example}

More systematically, we have employed two upscaling rules $S,R$ and
two downscaling rules $D,Q$ which obey
	\begin{align*}
	SL=L^2 S,\quad RL=L^2 R,\quad DL^2=LD,\quad DQ^2=LQ
	\end{align*}
($L$ is left shift). We have data $p^{(j)}$ at
level $j$, $j=0,\dots,M$, where we interpret
the data at the highest (finest) level as given, and the data at lower
(coarser) level computed by downscaling. We also
store details $q^{(j)}$ at each level:
	\begin{align}
	p^{(j-1)} &= D p^{(j)}, 
	& q^{(j)} &= Q(p^{(j)}\ominus Sp^{(j-1)}).
	\label{jw:eq:decompose}
	\end{align}
We require that upscaled level $j-1$ data and level $j$ details can restore
level $j$ data:
	\begin{align}
	p^{(j)} = Sq^{(j-1)}\oplus Rq^{(j)}.
	\label{jw:eq:reconstruct}
	\end{align}
Generally, $S$, $D$ compute points from points, so they are formulated
via averages:
	\begin{align*}
	\Sp_i&=\avg(a_{i-2j};p_j), 
	&Dp_i&=\avg(a_{2i-j};p_j).
	\end{align*}
In Example~\ref{jw:ex:haar}, averages are computed w.r.t.\ base points
$p_{\lfloor i/2\rfloor}$ for $S$ resp.\ $p_i$ for $D$, and coefficients
$a_j$ and $b_j$ vanish except $a_0=a_1=1$, $b_{0}=b_{1}={1\over 2}$.

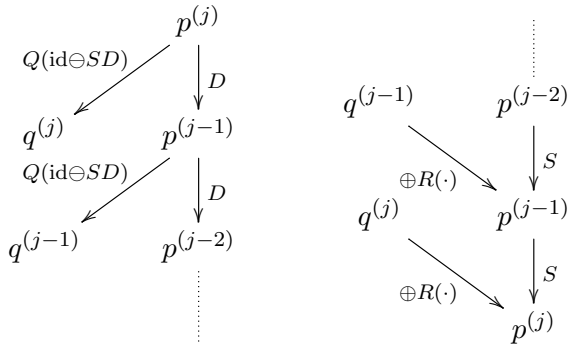
\begin{figure}[b]
\begin{minipage}{.25\textwidth}
\xymatrix{
		& p^{(j)} 
			\ar[ld]_{Q(\id\ominus SD)}
			\ar[d]^{D} \\
	q^{(j)}
		& p^{(j-1)}
			\ar[ld]_{Q(\id\ominus SD)}
			\ar[d]^{D} \\
	q^{(j-1)}
		& p^{(j-2)}
			\ar@{.}[d]
			\\
		&
}
\end{minipage}
\hfill
\begin{minipage}{.25\textwidth}
\xymatrix{
		&	
			\ar@{.}[d]
			\\
	q^{(j-1)}
	\ar[dr]_{\oplus R(\cdot)}
		& p^{(j-2)}
		\ar[d]^S
			\\
	q^{(j)}
	\ar[dr]_{\oplus R(\cdot)}
		& p^{(j-1)}
		\ar[d]^S
			\\
		& p^{(j)}
}
\end{minipage}
\hfill
	\begin{minipage}[b]{0.4\textwidth}
	\caption{The decomposition and reconstruction chains of operations
		in a geometric multiscale decomposition based on
		upscaling and downscaling $S,D$ for points, and 
		upscaling and downscaling $R,Q$ for
		detail vectors.}
	\end{minipage}
\end{figure}

The downscaling operator $Q$ acts on tangent vectors
$v_i = p_i\ominus (SDp)_i\in T_{p_i}\MM$, so it has to deal with
vectors potentially contained in different vector spaces. 
In our special case, $Q$ simply forgets one half of the data:
	\begin{align*}
	(Qp)_i = p_{2i}.
	\end{align*}
Finally, the upscaling operator $R$ takes the vectors stored in $q^{(j)}$
and converts them into vectors which can be added to upscaled points
$Sp^{(j-1)}$. Thus $R$ potentially has to deal with vectors contained
in different tangent spaces. In our special case, the points
$(Sp^{(j-1)})_{2i}$ $(Sp^{(j-1)})_{2i+1}$ both coincide with $p^{(j-1)}_i$,
and that is also the point where the detail coeffient $q^{(j)}_i$ is
attached to. We therefore might be tempted to write
$(Rq)_{2i} = q_i$, $(Rq)_{2i+1} = -q_i$. This simple rule
however does not take into account that along reconstruction, 
data and details could have been modified, and no longer fit together.
We therefore use parallel transport to move the vector to the right
tangent space, just in case:
	\begin{align*}
	(Rq)_{2i} &= \PT^{(p^{(j-1)})_i} (q_i),
	&(Rq)_{2i+1} &= - (Rq)_{2i}.
	\end{align*}
 The symbol $\PT^{(p^{(j-1)})_i}(q_i)$ refers to transporting $q_i$ to
a tangent vector attached to $({p^{(j-1)}})_i$, see
Section~\ref{jw:sec:hermite}.

The operations $S,R,D,Q$ must be compatible, in the sense that
reconstruction is a left inverse of downscaling plus computing details.
While in the linear case, where $S,D,R,Q$ are linear operators on
$\ell^\infty(\RR^d)$, one usually requires $QR=\id$ and $QS=0$ as well
as $SD+RQ=\id$, in the geometric case we must be careful not to mix
operations on points with operations on tangent vectors. We therefore
require 
	\begin{align}
	SDp \oplus (RQ(p\ominus SDp)) = p.
	\label{jw:eq:compat}
	\end{align}

\index{wavelet!interpolatory}
\begin{Example}{Interpolatory wavelets}\label{jw:ex:interp:wav}\relax
Consider an interpolatory subdivision rule $S$ with dilation factor $2$,
i.e., $\Sp_{2i}=p_i$, and the forgetful downscaling operator 
$p^{(j-1)}_i = (Dp^{(j)})_i=p^{(j)}_{2i}$.
If we store as details the difference vectors
between $SDp$ and $p$ for odd indices, the data points $p_{2i+1}$ can
be easily reconstructed:
	\begin{align*}
	p^{(j-1)}_i &= p^{(j)}_{2i},\
	&q^{(j)}_i &= p^{(j)}_{2i+1} \ominus Sp^{(j-1)} 
			& \text{(decomposition)}, \\
	p^{(j)}_{2i} &= p^{(j-1)}_i, \
	&p^{(j)}_{2i+1} &= (Sp^{(j-1)})_{2i+1}\oplus q^{(j)}_i
			& \text{(reconstruction)}.
	\end{align*}
 This procedure fits into the general scheme described above if we
we let $Q=DL$ ($L$ is left shift) and define the upscaling of
details by $(Rq)_{2i}=0$, $(Rq)_{2i+1}=q_i$. To admit the
possibility that before reconstruction, data and details have
been changed, we define
	\begin{align*}
	\index{parallel transport}
	(Rq)_{2i}&=0 \in T_x\MM, 
	&(Rq)_{2i+1} &= \PT^x(q^{(j)}_i)\in T_x\MM,
	&\text{where}\ x&= {Sp^{(j-1)}_{2i+1}},
	\end{align*}
in order to account for the possibility that $q_i^{(j)}$ is not
yet contained in the ``correct'' tangent space. 
The decimated data $p^{(j-1)}$ together with details $q^{(j)}$
($j\le M$) may be called a geometric interpolatory-wavelet decomposition
of the data at the finest level $p^{(M)}$. That data itself comes e.g.\
from sampling a function, cf.\ \cite{jw:donoho-tr}. 
\end{Example}

\Paragraph{Definability of multiscale transforms without redundancies}
The previous examples use upscaling and downscaling operations which
are rather simple, except that in Example~\ref{jw:ex:interp:wav}
one may use any interpolatory subdivision rule. It is also possible
to extend Example \ref{jw:ex:haar} to the more general case of 
a midpoint-interpolating subdivision rule $S$,
which is a right inverse of the decimation operator $D$.
In \cite{jw:grohs-2012-ds} it is argued that it is highly
unlikely that in the setup described above, which avoids
redundancies, more general upscaling and downscaling rules will manage
to meet the compatibility condition \eqref{jw:eq:compat} needed for
perfect reconstruction. In the linear case, where all details
$q_i^{(j)}$ are contained in the {\em same} vector space, 
\eqref{jw:eq:compat} is merely an algebraic condition on 
the coefficents involved in the definition of $S,D,Q,R$ which can be solved.
In the geometric case, the usage of parallel transport makes a fundamental
difference.

\subsection{Properties of multiscale transforms}

\Paragraph{Charcterizing smoothness by coefficient decay}
One purpose of a multiscale decomposition of data is to read off 
properties of the original data. Classically, the faster the
magnitude of detail coefficients $q^{(j)}_i$ decays as $j\to\infty$,
the smoother the original data. A corresponding result for
the interpolatory wavelets of Example~\ref{jw:ex:interp:wav} 
in the linear case is given by \cite[Th. 2.7]{jw:donoho-tr}. 
To state a result in the multivariate geometric case, let us first
introduce new notation for interpolatory wavelets,
superseding Example~\ref{jw:ex:interp:wav}.

We consider an interpolatory subdivision rule $S$ acting
with dilation factor $N$ on data $p\colon\ZZ^s\to\MM$. 
We define data $p^{(j)}$ at level $j$ as samples of a
function $f\colon\RR^s\to\MM$, and construct detail vectors 
similar to Example~\ref{jw:ex:interp:wav}:
	\begin{align}
	p^{(j)}_i &= f(N^{-j} i),
	& q^{(j)} & = p^{(j)}\ominus  Sp^{(j-1)},
	& p^{(j)} &= Sp^{(j-1)} \oplus q^{(j)}.
	\label{jw:eq:intp:wav2}
	\end{align}
 This choice is consistent with the decimation operator $Dp_i:=p_{Ni}$.
The difference to Example~\ref{jw:ex:interp:wav} is firstly that
here we allow multivariate data, and secondly that we do not ``forget'' 
redundant information such as $q^{(j)}_{Ni}=0$. 

The result below uses  the notation
$\Lip\gamma$ for functions which are $C^k$ with
$k=\lfloor\gamma\rfloor$ and whose
$k$-th derivatives are H\"older continuous of exponent $\gamma-k$.
The critical H\"older regularity of a function $f$ is the supremum
of $\gamma$ such that $f\in\Lip\gamma$.

\begin{Satz}{Theorem} {\rm \cite[Th. 8]{jw:grohs-2009-wav}}
Assume that the interpolatory upscaling rule $S$,
when acting linearly on data $p\colon\ZZ^s\to\RR$,
reproduces polynomials of degree $\le d$ and has limits of
critial H\"older regularity $r$.

Consider a continous function $f\colon\RR^s\to\MM$, and construct
detail vectors $q^{(j)}$ at level $j$ for the function
$x\mapsto f(\sigma\cdot s)$ for some $\sigma>0$ (whose local existence
is guaranteed for some $\sigma>0$).

Then $f\in\Lip\alpha$, $\alpha<d$ implies that detail vectors decay
with $\sup_i\|q^{(j)}_i\| \le C\cdot N^{-\alpha j}$ as $j\to\infty$.
Conversely, that
decay rate together with $\alpha<r$ implies $f\in \Lip\alpha$.  The
constant is understood to be uniform in a compact set.
\end{Satz}

The manifold $\MM$ can be any of the cases we defined $\oplus$ and
$\ominus$ operations for. Of course, smoothness of $f\colon\RR^s\to\MM$
is only defined up to the intrinsic smoothness of $\MM$ as
a differentiable manifold.
An example of an upscaling rule $S$
is the four-point scheme with parameter $1/16$  mentioned in
Example~\ref{jw:ex:conv:interpolatory}, which
reproduces cubic polynomials and has critical H\"older regularity 2, cf.\
\cite{{jw:dyn:1987:4p}}.

The proof is conducted in a coordinate chart (it does not matter which),
and uses a linear vision of the theorem as an auxiliary tool.
It further deals with the 
extensive technicalities which surround proximity
inequalities in the multivariate case.

It is worth noting that A.\ Weinmann 
in \cite{jw:weinmann-2012-mf} succeeded in transferring these ideas
to the combinatorially irregular setting. The results are essentially
the same, with the difference that one can find upscaling rules only
up to smoothness $2-\epsilon$.

\Paragraph{Stability}
Compression of data is a main application of multiscale decompositions,
and it is achieved e.g.\ by thresholding or quantizing detail vectors.
It is therefore important to know what effect these changes have 
when reconstruction is performed. What we bascially want to know is
whether reconstruction is Lipschitz continuous. In the linear case
the problem does not arise separately, since the answer is implicitly
contained in norms of linear operators. For the geometric
multiscale transforms
defined by upscaling operations $S,R$ and downscaling operations $D,Q$
according to \eqref{jw:eq:compat}, this problem is discussed by 
\cite{jw:grohs-2012-ds}. 
Consider data $p^{(j)}$ at level
$j$ with $p^{(j-1)}=Dp^{(j)}$ such that $\delta(p^{(j)}) \le
C\mu^j$, for some $\mu<1$. Consider recursive reconstruction of data
$p^{(j)}$ from $p^{(0)}$ and details $q^{(1)},\ldots,q^{(j)}$ 
according to Equation \eqref{jw:eq:reconstruct}. Then there are
constants $C_k$ such that for modified details $\tilde q^{(j)}$,
leading to modified data $\tilde p^{(j)}$, we have the local
Lipschitz-style estimate 
	\begin{align*}
	& \sup\nolimits_i\|p^{(0)}_i-\tilde p^{(0)}_i\| \le C_1,
	\sup\nolimits_i\|q^{(k)}_i-\tilde q^{(k)}_i\| \le C_2\mu^k
	\\
	\implies &
	\sup\nolimits_i\|p^{(j)}_i-\tilde p^{(j)}_i\| \le
	C_3\Big(
		\sup\nolimits_i
			\|p^{(0)}_i-\tilde p^{(0)}_i\|
		+\sum\nolimits_{k=1}^j 
			\sup\nolimits_i
			 \|q^{(k)}_i-\tilde q^{(k)}_i\|
		\Big).
	\end{align*}
It refers to a coordinate chart of
the manifold $\MM$ (it does not matter which).

\index{stability}

\Paragraph{Approximation Order}
For an interpolatory upscaling operator $S$, and 
data $p_i\in\MM$ defined by sampling,
$p_i = f(h\cdot i)$, we wish to know to 
what extent the original function differs
from the limit created by upscaling the sample. 
We say that $S$ has approximation order $r$, if
there are $C>0$, $h_0>$ such that for all $h<h_0$ 
	\begin{align*}
	\sup\nolimits_x\dMM(S^\infty f(x/h), f(x))
	\le C\cdot h^r. 
	\index{approximation order}
	\end{align*}
It was shown by
\cite{jw:yu:approximationorder-2011} that a manifold subdivision
rule has in general the same approximation order as the linear
rule we get by restricting $S$ to linear data.

This question is directly related to stability as 
discussed above:
Both $f$ and $S^\infty p$ can be reconstructed 
from samples $p^{(j)}$, if $h=N^{-j}$: Detail vectors $q^{(k)}$, $k>j$,
according to \eqref{jw:eq:intp:wav2} reconstruct $f$, whereas details
$\tilde q^{(k)}=0$ reconstruct $S^\infty p$.
Stability of reconstruction and knowledge of the
asymptotic magnitude of details $q^{(k)}_i$, $k>j$ directly corresponds
to approximation order. On basis of this relationship 
one can again show an approximation order equivalence result, cf.\
\cite{jw:grohs-2010-ao}.

\subsubsection*{Conclusion}
The preceding pages give an account of averages, subdivision, and
multiscale transforms defined via geometric
operations which are intrinsic for various geometries (metric spaces,
Riemannian manifolds, Lie groups, and symmetric spaces). We 
reported on complete solutions 
in special cases (e.g.\ convergence of subdivision rules
in Hadmard metric spaces) and on other results with much more general
scope as regards the spaces and subdivision rules involved, but with
more restrictions on the data they apply to.

\subsection*{Acknowledgments}

The author gratefully acknowledges the support of the Austrian
Science fund through grant No.\ W1230, and he also wants to thank
all colleagues whose work he was able to present in this survey.

\index{subdivision|)}
\index{multiscale transform|)}


\providecommand{\bysame}{\leavevmode\hbox to3em{\hrulefill}\thinspace}
\providecommand{\MR}{\relax\ifhmode\unskip\space\fi MR }
\providecommand{\MRhref}[2]{%
  \href{http://www.ams.org/mathscinet-getitem?mr=#1}{#2}
}
\providecommand{\href}[2]{#2}

\end{document}